\documentclass[12pt,a4paper]{article}
\usepackage[english]{babel}
\usepackage[psamsfonts]{amssymb}
\usepackage{cmmib57}
\usepackage{exscale}
\usepackage{amsmath}
\usepackage{amscd}
\usepackage{latexsym}
\usepackage{graphicx}
\usepackage{colortbl}
\usepackage{pdfsync}
\numberwithin{equation}{section}
\numberwithin{subsection}{section}
\allowdisplaybreaks
\numberwithin{equation}{section}
\newtheorem{stat}{Statement}[section]
\newtheorem{theorem}[stat]{Theorem}
\newtheorem{corollary}[stat]{Corollary}
\newtheorem{proposition}[stat]{Proposition}
\newtheorem{lemma}[stat]{Lemma}

\newcommand{\cqd}{\hfill$\Box$}

\newcommand{\Rd}{{\mathbb{R}^d}}
\newcommand{\R}{{\mathbb R}}
\newcommand{\N}{{\mathbb N}}
\newcommand{\D}{{\mathbb D}}
\newcommand{\tf}{\mathcal{F}}
\newcommand{\caA}{\mathcal{A}}

\newcommand{\beq}{\begin{equation}}
\newcommand{\eeq}{\end{equation}}
\newcommand{\bal}{\begin{align}}
\newcommand{\eal}{\end{align}}
\newcommand{\beqn}{\begin{equation*}}
\newcommand{\eeqn}{\end{equation*}}
\newcommand{\baln}{\begin{align*}}
\newcommand{\ealn}{\end{align*}}

\newcommand{\ep}{\epsilon}

\newcommand{\Pb}{\mathbb P}
\newcommand{\E}{\mathbb{E}}
\newcommand{\hac}{\mathcal{H}}

\newcommand{\tx}{{(t,x)\in[0,T]\times \mathbb{R}^d}}

\begin{document}
%%%%%%%%%%%%%%%
\begin{titlepage}
\vskip 1cm

\begin{center}
{\Large\bf The stochastic wave equation in high dimensions: Malliavin differentiability and absolute continuity}
\medskip

by\\
\vspace{14mm}

\begin{tabular}{l@{\hspace{10mm}}l@{\hspace{10mm}}l}
{\sc Marta Sanz-Sol\'e}$\,^{(\ast)}$ &and &{\sc Andr\'e S{\"u}{\ss}}$\,^{(\ast)}$\\
{\small marta.sanz@ub.edu}         &&{\small andre.suess@ub.edu}\\
{\small http://www.mat.ub.es/$\sim$sanz}\\
\end{tabular}
\begin{center}
{\small Facultat de Matem\`atiques}\\
{\small Universitat de Barcelona } \\
{\small Gran Via de les Corts Catalanes, 585} \\
{\small E-08007 Barcelona, Spain} \\
\end{center}

\vskip 0.1cm
\today.

\end{center}

\vskip 1cm

\noindent{\bf Abstract} We consider the class of non-linear stochastic partial differential equations studied in \cite{conusdalang}. Equivalent formulations using integration with respect to a cylindrical Brownian motion and also the Skorohod integral are established. It is proved that the random field solution to these equations at any fixed point $(t,x)\in[0,T]\times \Rd$ is differentiable in the Malliavin sense. For this, an extension of the integration theory in \cite{conusdalang} to Hilbert space valued integrands is developed, and commutation formulae of the Malliavin derivative and stochastic and pathwise integrals are proved. In the particular case of equations with additive noise, we establish the existence of density for the law of the solution at $(t,x)\in]0,T]\times\Rd$. The results apply to the stochastic wave equation in spatial dimension $d\ge 4$.

  \vskip 1cm

\noindent{\bf Keywords:} stochastic integration; stochastic partial differential equations; stochastic wave equation; Malliavin calculus; densities.

\smallskip

\noindent{\bf MSC 2000:} Primary: 60H15, 60H07; Secondary: 60H20, 60H05. 
\vfill

\noindent
\footnotesize
{\begin{itemize} \item[$^{(\ast)}$] Supported by the grant MICINN-FEDER MTM 2009-07203 from the \textit{Direcci\'on General de
Investigaci\'on, Ministerio de Economia y Competitividad, Spain.}
\end{itemize}}

\end{titlepage}
\newpage

%\end{document}

%%%%%%%%%%%
%%%%%%%INTRODUCTION
\section{Introduction}
\label{s0}
In this article, we consider stochastic partial differential equations (SPDEs) of the type
\begin{align}
\label{i.0}
L u(t,x) &= \sigma(u(t,x)) \dot F(t,x) + b(u(t,x)),\ (t,x)\in]0,T]\times \Rd,\nonumber\\
u(0,x) &= \frac{\partial}{\partial t}u(0,x)=0, \ x\in\Rd,
\end{align}
where $L$ is a  second order differential operator, $\sigma$ and $b$ are real functions, $\dot F$ is the formal derivative
of a Gaussian stochastic process and $d\in \N$. The setting applies in particular to the wave operator
\begin{equation*}
\frac{\partial^2}{\partial t^2}- \Delta_d,
\end{equation*}
where $\Delta_d$ denotes the Laplacian in dimension $d$. 

We give a rigorous meaning to \eqref{i.0} using the stochastic integration theory developed in \cite{conusdalang}, which extends the setting of \cite{walsh} and \cite{dalang}. More precisely, let $G$ be the fundamental solution associated with the operator $L$. We consider the mild form of the equation \eqref{i.0},
\begin{align}
\label{i.1}
u(t,x)&=\int_0^t \int_{\Rd} G (t-s,x-z) \sigma(u(s,z)) M(ds,dz)\nonumber\\
&+ \int_0^t \int_{\Rd} G (t-s,x-z) b(u(s,z)) dz ds,
\end{align}
where $M$ denotes the martingale measure derived from the Gaussian process $F$. More specifications on $F$ and on the stochastic and pathwise integrals used in \eqref{i.1} are given in Section \ref{s1}. Under suitable conditions on $G$, and for Lipschitz continuous coefficients $\sigma$ and $b$, \cite[Theorem 4.2]{conusdalang} establishes the existence of a random field solution to \eqref{i.1}. 

The main objective is to establish the differentiability in the Malliavin sense of the random variable $u(t,x)$ defined by \eqref{i.1}, for each fixed $(t,x)\in[0,T]\times\Rd$. As a consequence, we will obtain results on the existence of density for $u(t,x)$, $(t,x)\in]0,T]\times\Rd$ which are applied to a stochastic wave equation in spatial dimension $d\ge 4$.

Malliavin differentiability and existence and properties of the density have been studied in particular cases of \eqref{i.1}, like the stochastic heat equation with $d\in\N$ and the stochastic wave equation with $d\in\{1,2,3\}$. We refer the reader to \cite{ballypardoux}, \cite{mcms}, \cite{milletsanzsole}, \cite{nualartsanzsole1}, \cite{nualartsanzsole2}, \cite{nualartquer}, \cite{pardouxtusheng}, \cite{quersanz}, \cite{quersanz2}, \cite{sanz}, \cite{sanzbook} for a sample of results. However, to the best of our knowledge, similar problems for the stochastic wave equation
in dimension $d\ge 4$ have not been so far solved. The main difficulty stems from the non-smoothness of the distribution $G$, the fundamental solution associated with the differential operator $L$. The results of this paper are general enough to cover that important example.

Next, we describe the content of the article. Section \ref{s1} gathers the preliminary notions and results underpinning Equation \eqref{i.0}, following \cite{conusdalang}. In particular, the governing noise $F$ is described and the construction of the integrals in \eqref{i.1} sketched. Along with this, we prove a new result (see Lemma \ref{lem:P0=HT}) that will be used later on in Section \ref{s3} to give a formulation of the stochastic integral in \cite{conusdalang} in terms of a sequence of independent standard Brownian motions. This provides a better understanding of this integral and makes its handling easier.

In Section \ref{s2} an extension of Conus-Dalang's stochastic and pathwise integrals to Hilbert space valued stochastic processes is developed. With this, we extend the results proved in \cite{quersanz} and provide the theoretical background for the study of the Malliavin derivative of the solution of \eqref{i.1}. 

Section \ref{s3} contains some complements to the Conus-Dalang's stochastic integral. For a relevant class of integrands, we prove that its divergence operator (in Malliavin sense) coincides with that integral and also with an It\^o stochastic integral with
respect to a sequence of independent standard Brownian motions. Actually, the latter is nothing but the stochastic integral with respect to a cylindrical Brownian motion, as in the setting of \cite{dapratozabczyk}. In particular, a partial extension of  \cite[Proposition 2.6]{dalangquer} is obtained.  

Section \ref{s4} contains preliminaries to Section \ref{s5}. It is proved that under suitable hypotheses, the Malliavin operator $D$ commutes with the stochastic integrals of \cite{conusdalang}. Then, in Section \ref{s5}, we prove that for any $(t,x)\in[0,T]\times \Rd$, the random variable $u(t,x)$ defined in \eqref{i.1} belongs to $\mathbb{D}^{1,2}$ - the space of Malliavin differentiable random variables with square integrable
derivative. Using a standard approach, we consider a sequence of $L^2(\Omega)$-approximations of the process $\{u(t,x), (t,x)\in[0,T]\times \Rd\}$, $\{u_n(t,x), (t,x)\in[0,T]\times \Rd\}$, $n\in\N$, such that they are Malliavin differentiable with uniformly bounded (in the corresponding norm) Malliavin derivatives. 
The section ends by establishing an SPDE satisfied by the Hilbert space valued process $Du(t,x)$ (see \eqref{5.1}). 
In contrast with examples where the distribution $G$ is smooth (for example, either a function, as for the stochastic heat equation, or a positive measure, as for the stochastic wave equation in dimension $d\in\{1,2,3\}$), it seems not possible to obtain this equation by passing to the limit the sequence of SPDEs \eqref{5.3} satisfied by the Malliavin derivatives of the approximations $u_n$. Indeed, for this one would need to have $u(t,x)\in\mathbb{D}^{1,p}$, for some $p>2$, a property that has not been established yet. This problem stems from the lack of $L^p(\Omega)$ estimates for the solution of \eqref{i.1} pointed out in \cite{conusdalang}. So far, this has been only proved when
$\sigma$ is an affine function (see \cite[Section 6]{conusdalang}). The general case is by now an open problem.

We overcome this problem and eventually establish \eqref{5.1}, by applying the operator $D$ to Equation \eqref{i.1} and the commutation results of Section \ref{s4}.
The stochastic integrals  in \cite{dalang} and \cite{conusdalang} are constructed assuming the property of spatial stationary covariance of the integrand. For the stochastic integral in \eqref{i.1} this follows from the ``$S$" property introduced in \cite{dalang} (see \cite[Definition 4.4 and Lemma 4.5]{conusdalang}). In the application of the results proved in Section \ref{s4}, we need in addition the stationary covariance property for Hilbert space valued stochastic processes of the form $\{D[B(u(t,x))], (t,x)\in[0,T]\times \Rd\}$, where $B$ is a smooth function. This can be achieved by considering 
the equation \eqref{5.a}, more general  than \eqref{i.1},
 and by proving that the ``$S$" property holds for the couple consisting of the solutions to these two equations. %This is to obtain the crucial result proved in Lemma \ref{5.3}.
 
The final Section \ref{s6} deals with the existence of density for each random variable $u(t,x)$, $(t,x)\in]0,T]\times \Rd$ in \eqref{i.1} when $\sigma$ is constant. The results apply to the stochastic wave equation in spatial dimension $d\ge4$ with an additive, Gaussian, spatially correlated  noise, with covariance measure given by a Riesz kernel. This is proved by applying the Bouleau-Hirsch's criterion. The reason for the restriction to additive noise comes from the fact that so far we have not been able to obtain lower bounds of the dominant term of the Malliavin matrix for non-constant coefficients $\sigma$. Among the difficulties we encounter to solve this problem are the lack of $L^p$ estimates we alluded before and the lack of positivity of $G$.

Throughout this article, we shall use the usual convention of calling constants by the same letter, although they may vary from one expression to the other.

%%%%%%%%%% End Introduction

%%%%%%%%%%
%%%%%%%%%% SECTION PRELIMINARIES

\section{Preliminaries}
\label{s1}

In this section, we fix some notation, we present some general facts and we recall the stochastic and pathwise integrals from \cite{conusdalang} that will be used throughout the paper. 
The relevant spaces are described and some relationships between them are proved. 
 
Denote by  $\mathcal {C}^\infty_0(\Rd)$  the space of infinitely differentiable functions with compact support; $\mathcal{S}(\Rd)$ will denote the Schwartz space of rapidly decreasing $\mathcal{C}^\infty$ functions, $\mathcal{S}^\prime(\Rd)$ the space of tempered Schwartz distributions
and $\mathcal{S}^\prime_r(\Rd)$ the space of tempered Schwartz distributions with rapid decrease (see \cite{schwartz}).

Let $\zeta\in \mathcal{C}^\infty_0(\Rd)$ be nonnegative, with support included in the unit ball of $\Rd$ satisfying $\int_{\Rd}\zeta(x) dx=1$. Set $\zeta_n(x) := n^d\zeta(nx)$, $n\in\N$. Then, 
as $n\to\infty$, $\zeta_n\rightarrow\delta_0$ in ${\mathcal{S}}^\prime(\Rd)$, and $\tf\zeta_n\rightarrow 1$ pointwise. Moreover $|\tf\zeta_n| \leq 1$ for all $n\in\N$. The sequence $(\zeta_n)_{n\in \N}$ is termed an {\it approximation of the identity}.

Let $\Lambda\in\mathcal{S}_r^\prime(\Rd)$. Denote by ``$\ast$"  the convolution operation. It is well-known that 
\begin{equation}
	\Lambda_n:=\Lambda\ast\zeta_n,
\label{eq:Lambda_n}
\end{equation}	
 belongs to $\mathcal{S}(\Rd)$. Moreover,
\begin{equation}
 \left\vert\tf\Lambda_n(\xi)\right\vert = \left\vert \tf\Lambda(\xi)\right\vert \left\vert\tf\zeta_n(\xi)\right\vert \leq \left\vert\tf\Lambda(\xi)\right\vert,
\label{eq:fourierpsi}
 \end{equation}
for all $\xi\in\Rd$, and  $\tf\Lambda_n\rightarrow\tf\Lambda$ pointwise as $n\to\infty$.

%%%%% Description of the noise
Let $\{F(\phi);\; \phi\in\mathcal{C}_0^\infty(\mathbb{R}_+\times\Rd)\}$ be a Gaussian process with mean zero and covariance functional
\begin{equation}
	\E[F(\phi)F(\psi)] = \int_0^\infty\int_\Rd \big(\phi(t)\ast\tilde{\psi}(t)\big)(x)\Gamma(dx) dt,
	\label{eq:correlation}
\end{equation}
where $\tilde{\psi}(t,x) := \psi(t,-x)$ and $\Gamma$ is a nonnegative, nonnegative definite, tempered measure on $\Rd$. There exists a nonnegative tempered measure $\mu$ on $\Rd$ such that $\tf\mu = \Gamma$ (see for instance \cite{schwartz}, Chapter VII, Th\'{e}or\`{e}me XVIII). Then by Parseval's identity, the right-hand side of \eqref{eq:correlation} is equal to  
\begin{equation*}
\int_0^{\infty}\int_{\Rd}\tf \phi(t)(\xi)\overline{\tf \psi(t)}(\xi)\mu(d\xi) dt.
\end{equation*}

 As is explained in \cite{dalangfrangos}, the process $F$ can be extended to a worthy martingale measure 
 $M=(M_t(A);\; t\in\mathbb{R}_+, A\in\mathcal{B}_b(\Rd))$ where $\mathcal{B}_b(\Rd)$ denotes the bounded Borel subsets of $\Rd$. The natural filtration generated by this martingale measure will be denoted in the sequel by $(\mathcal{F}_t)_{t\ge 0}$.
 
  A stochastic integration theory with respect to martingale measures has been developed by M\'etivier and Pellaumail and by Walsh, among others. Here, we shall use \cite{walsh} as reference. Using this integral, we have
 \begin{equation}
	F(\phi) = \int_0^T\int_\Rd \phi(s,z)M(ds,dz),
	\label{eq:FandM}
\end{equation}
for $\phi\in\mathcal{C}_0^\infty(\mathbb{R}_+\times\Rd)$.
 
 Extensions of the stochastic integral given in \cite{walsh} have been introduced in \cite{dalang} and more recently, in \cite{conusdalang}. Throughout this article, we shall refer mainly to the latter.

%Integrands

Fix $T>0$. For stochastic processes $f$ and $g$, indexed by $(t,x)\in[0,T]\times\Rd$ and satisfying suitable conditions, we define the inner product
\begin{align*}
 \langle f,g\rangle_{0}& = \E\left[\int_0^T\int_\Rd \big(f(s)\ast\tilde{g}(s)\big)(x)\Gamma(dx) ds \right]\\
 & = \E\left[\int_0^T\int_\Rd \tf f(s)(\xi)\overline{\tf g(s)(\xi)}\mu(d\xi) ds\right], 
 \end{align*}
where the corresponding norm $\|\cdot\|_0$ is defined in the usual way. Moreover, we define the norm
\[ \|g\|_{+}^2 = \E\left[\int_0^T\int_\Rd \big(|g(s)|\ast|\tilde{g}(s)|\big)(x)\Gamma(dx) ds\right]. \]
Let $\mathcal{P}_+$ be the set of predictable processes $g$ such that  $\|g\|_+ < \infty$. In \cite[Exercise 2.5]{walsh} it is shown that $\mathcal{P}_+$ is complete and hence it is a Banach space. 
Let $\mathcal{E}$ denote the set of {\it simple processes} $g$, that is, stochastic processes of the form 
\begin{equation}
	g(t,x;\omega) = \sum_{j=1}^m 1_{(a_j,b_j]}(t)1_{A_j}(x)X_j(\omega),
	\label{eq:simpleprocess}
\end{equation}  
for some $m\in\N$, where $0\leq a_j < b_j\leq T$, $A_j\in\mathcal{B}_b(\Rd)$ and $X_j$ is a bounded and $\mathcal{F}_{a_j}$-measurable random variable for all $1\leq j\leq n$.  

%\[ (g\cdot M)_t(B) := \sum_{j=1}^n X_j\big(M_{b\wedge t}(A_j\mathcal{P} B) - M_{a\wedge t}(A_j\mathcal{P} B)\big), \]
According to \cite[Proposition 2.3]{walsh}, $\mathcal{E}$ is dense in $\mathcal{P}_+$. Hence, we can also define $\mathcal{P}_+$ as the completion of $\mathcal{E}$ with respect to $\|\cdot\|_+$.

Following \cite{dalang}, we denote by $\mathcal{P}_0$ the completion of $\mathcal{E}$ with respect to $\|\cdot\|_0$. This is a Hilbert space consisting of predictable processes which contains tempered distributions in the $x$-argument (whose Fourier transform are functions, $\Pb$-a.s.).  The norm in this space is given by
\begin{equation}
	\|g\|^2_0 = \E\bigg[\int_0^T\int_\Rd |\tf g(s)(\xi)|^2 \mu(d\xi)ds\bigg].
	\label{eq:normP01}
\end{equation}
For sufficiently smooth elements of $\mathcal{P}_0$, this norm can be also written as
\begin{equation}
	\|g\|^2_0 = \E\bigg[\int_0^T\int_\Rd \big(g(s,\cdot)\ast\tilde{g}(s,\cdot)\big)(z)\Gamma(dz)ds\bigg].
	\label{eq:normP02}
\end{equation}
Note that $\mathcal{P}_0$ is not defined as the set of predictable processes for which $\|\cdot\|_0 < \infty$. In fact, it can be shown that the latter space is not complete. Since $\|\cdot\|_0 \leq \|\cdot\|_+$, we clearly have $\mathcal{P}_+\subseteq\mathcal{P}_0$, and from the above comments on completeness we know that this inclusion must be strict.
\smallskip

Consider the subsets of $\mathcal{P}_+$  and $\mathcal{P}_0$ consisting of deterministic processes, denoted by 
$\mathcal{P}_{+,d}$ and $\mathcal{P}_{0,d}$, respectively. In the next lemma, we give an equivalent definition of $\mathcal{P}_{0,d}$. For this, we first introduce a new space $\mathcal{E}_0$ consisting of Schwartz functions endowed with the inner product 
\begin{equation}
 \langle \phi,\psi\rangle_0 = \int_\Rd \big(\phi\ast\tilde{\psi}\big)(x)\Gamma(dx) = \int_\Rd  \tf\phi(\xi)\overline{\tf\psi(\xi)}\mu(d\xi),
\label{eq:spHT}
\end{equation}
where $\phi,\psi\in\mathcal{S}(\Rd)$. Let $\mathcal{H}$ denote the completion of $(\mathcal{E}_0,\langle \cdot,\cdot\rangle_0)$
and set ${\mathcal{H}_T} := L^2([0,T];\mathcal{H})$. In the sequel we will also denote by $\Vert\cdot\Vert_{\mathcal{H}_T}$ the norm in this space derived from the scalar product $\langle\cdot,\cdot\rangle_0$.
%%%%
%%%%
%%%% Relations between spaces

\begin{lemma}\label{lem:P0=HT}
The spaces $\mathcal{P}_{0,d}$ and $\mathcal{H}_T$ coincide. 
\end{lemma}
\noindent{\it Proof.} First, we prove the inclusion $\mathcal{H}_T\subseteq \mathcal{P}_{0,d}$. For this, 
let $\mathcal{E}_s$ be the set of functions
$\phi:[0,T]\times\Rd\rightarrow\R$ which are a step function in the first argument and a Schwartz function in the second one. Notice that $\mathcal{E}_s$
is dense in ${\mathcal{H}_T}$ with respect to the norm $\Vert\cdot\Vert_{\mathcal{H}_T}$.
We will show  that $\mathcal{E}_s\subseteq\mathcal{P}_{+,d}\subseteq\mathcal{P}_{0,d}$, yielding the statement. 

Indeed, fix $\phi\in\mathcal{E}_s$. Due to Leibniz' formula (see \cite[Exercise 26.4]{treves}), the function $z\mapsto\big(|\phi(s,\cdot)|\ast|\tilde{\phi}(s,\cdot)|\big)(z)$
decreases faster than any polynomial in $|z|^{-1}$. Since $\Gamma$ is a tempered measure, we have
\[ \|\phi\|_+^2 =  \int_0^T\int_\Rd \big(|\phi(s,\cdot)|\ast|\tilde{\phi}(s,\cdot)|\big)(z)\Gamma(dz)ds < \infty. \]
This proves the claim.

Next, we consider the set $\mathcal{E}_d$ consisting of deterministic simple functions, and we prove that
 $\mathcal{E}_d\subseteq\mathcal{H}_T$. By taking closures in the norm $\Vert\cdot\Vert_{\mathcal{H}_T}$, we will obtain the inclusion
 $\mathcal{P}_{0,d}\subseteq\mathcal{H}_T$.
 
Let $\psi\in\mathcal{E}_d$ be given by $\psi=1_{(a,b]}1_{A}$, with $0\leq a < b \leq T$. This function satisfies 
\begin{equation}
\label{1}
\Vert \psi\Vert^2_{\mathcal{H}_T}=\int_0^T \int_{\Rd}\left\vert \tf \psi(s)(\xi)\right\vert^2 \mu(d\xi) ds <\infty.
\end{equation}
Indeed, by writing $\Vert \psi\Vert^2_{\mathcal{H}_T}$ as in the right-hand side of \eqref{eq:normP02}, we have
\begin{align*}
\Vert \psi\Vert^2_{\mathcal{H}_T}&= \int_0^T\int_\Rd  \big((1_{(a,b]}(s)1_{A}(\cdot))\ast(1_{(a,b]}(s)\widetilde{1_{A}}(\cdot))\big)(z) \Gamma(dz)ds \\
  & = \int_0^T 1_{(a,b]}(s)ds\int_\Rd\int_\Rd 1_{A}(y)1_{A}(y-z)  dy\Gamma(dz)\\
  &\le (b-a)|A|\int_\Rd 1_{B}(z)\Gamma(dz) \leq (b-a)|A|\Gamma(\bar{B}),
  \end{align*}
  where $|A|$ denotes the Lebesgue measure of $A$, $B$ stands for  the ball in $\Rd$ centered at 0 and with radius ${\text{diam}}(A)=\sup\{d(x,z);\; x,y\in A\}$, and $\bar B$ denotes its closure in the Euclidean norm. Since $\Gamma$ is a nonnegative tempered measure, it has the form $\Gamma(dz) = p(z)\nu(dz)$, where $p$ is a polynomial and $\nu$ is a finite measure (see \cite[p.\ 242]{schwartz}). Hence $\Gamma$ is $\sigma$-finite. This fact along with the preceding inequalities, yields \eqref{1}.

For an approximation of the identity $(\zeta_n)_{n\in\N}$, we define $\psi_n(s):=\psi(s)\ast\zeta_n\in\mathcal{S}(\Rd)$, $s\in[0,T]$, $n\in \N$. Clearly, $\psi_n\in\mathcal{E}_0\subseteq{\mathcal{H}_T}$. 
Moreover, we will prove that
\begin{equation}
\label{2}
\lim_{n\to\infty}\Vert \psi_n-\psi\Vert_{\mathcal{H}_T}=0.
\end{equation}
This yields $\psi\in\mathcal{H}_T$.

For the proof of \eqref{2}, we notice that by the very definition of the norm in $\mathcal{H}_T$,
\begin{align}
  \|\psi_n-\psi\|_{\mathcal{H}_T}^2	& = \int_0^T\int_\Rd |\tf \psi_n(s)(\xi)-\tf \psi(s)(\xi)|^2 \mu(d\xi)ds \notag\\
                   		& = \int_0^T\int_\Rd |\tf \psi(s)(\xi)|^2|\tf\zeta_n(\xi)-1|^2 \mu(d\xi)ds  \label{eq:prfP0HT}.
                   		%& \leq 4 \int_0^T\int_\Rd |\tf g(s)(\xi)|^2|\mu(d\xi)ds.\ mathcal{H}_T
\end{align}

By using \eqref{1} and applying bounded convergence, the last term converges to zero as $n\to\infty$. 
\cqd
\medskip

Adding the random component yields the following.

\begin{corollary}\label{cor:P0=L2}
	The spaces $\mathcal{P}_0$ and the space of all predictable stochastic process in $L^2(\Omega\times[0,T];\hac)$ coincide.
\end{corollary}

%%%%%
%%%%%%%%%% reminders on stochastic integrals
%%%%%

In order to introduce notation and provide some introductory material, we give a brief overview of the integrals defined in \cite{conusdalang}. In the next section, we shall extend these integrals to Hilbert space valued stochastic processes. 

%%%%%%%Integrals
%%%% integrands

Let $Z=\{Z(t,x);\; t\in[0,T], x\in\Rd\}$ be a real-valued stochastic process, non identically zero, with the following properties. 
\begin{description}
\item{{\bf (A1)}}
$Z$ is a predictable stochastic process satisfying 
$\sup_\tx \E[Z(t,x)^2] < \infty$. 
\item{{\bf (A2)}}
$Z$ has  {\it spatial stationary covariance}. That is, for any $t\in[0,T]$, $x,y\in\Rd$, 
    \[ \E[Z(t,x)Z(t,x+y)] = \E[Z(t,0)Z(t,y)] =: \gamma^Z_t(y). \]
\end{description}

%In order to avoid technical subtleties when defining the norms in \eqref{eq:normP0Z1} below, we exclude the trivial case of $Z\equiv0$. 
The process 
\[ M_t^Z(A) := \int_0^t\int_A Z(s,z)M(ds,dz), t\in[0,T], A\in\mathcal{B}_b(\Rd),\]
defines a worthy martingale measure (see \cite{walsh}).

 Similarly to the definition of the norms $\|\cdot\|_+$ and $\|\cdot\|_0$, for stochastic processes indexed by $(t,x)\in[0,T]\times\Rd$ and satisfying suitable conditions, we set
\begin{align}
  \|g\|_{+,Z}^2 & = \E\bigg[\int_0^T\int_\Rd \big(|g(s,\cdot)Z(s,\cdot)|\ast|\tilde{g}(s,\cdot)\tilde{Z}(s,\cdot)|\big)(z)\Gamma(dz)ds\bigg], \notag\\
  \|g\|_{0,Z}^2 & = \E\bigg[\int_0^T\int_\Rd \big((g(s,\cdot)Z(s,\cdot))\ast(\tilde{g}(s,\cdot)\tilde{Z}(s,\cdot))\big)(z)\Gamma(dz)ds\bigg]. \label{eq:normP0Z1}
\end{align}
%where  the superscript ``\~\null" means that the spatial argument is changed into its opposite.
Let $\mathcal{P}_{+,Z}$ and $\mathcal{P}_{0,Z}$ denote the completion of $(\mathcal{E},\|\cdot\|_{+,Z})$ and $(\mathcal{E},\|\cdot\|_{0,Z})$ with respect to these norms, respectively. Accordingly to \cite[Exercise 2.5]{walsh} $\mathcal{P}_{+,Z}$ is exactly the set of all predictable processes $g$ for which $\|g\|_{+,Z}<\infty$. However, for $\mathcal{P}_{0,Z}$ there is no similar characterization. 

By Bochner's Theorem (\cite[Chapter VII, Th\'{e}or\`{e}me XVIII]{schwartz}), there exists a non-negative tempered measure $\nu_t^Z$ such that $\gamma^Z_t=\tf\nu_t^Z$, where 
$\gamma^Z_t$ is defined in {{\bf (A2)}}. Moreover, we have
 \begin{equation}
 \label{spectral}
 \gamma_t^Z\Gamma = (\tf\nu_t^Z) (\tf\mu) = \tf(\mu\ast\nu_t^Z). 
 \end{equation}

In the sequel we set $\mu_t^Z:=\mu\ast\nu_t^Z$. Due to Fubini's Theorem, assumption {{\bf (A2)}}, \eqref{spectral} and Parseval's Identity, we see that for all $g\in\mathcal{P}_{0,d}$
\begin{align}
  \|g\|_{0,Z}^2 
  %& = \E\bigg[\int_0^T\int_\Rd \int_\Rd g(s,y)Z(s,y)g(s,y-z)Z(s,y-z) dy\Gamma(dz)ds\bigg] \notag\\
  %& = \int_0^T\int_\Rd \int_\Rd g(s,y)g(s,y-z) dy\gamma^Z_s(z)\Gamma(dz)ds \notag\\
  & = \int_0^T\int_\Rd \big(g(s,\cdot)\ast \tilde{g}(s,\cdot)\big)(z)\gamma^Z_s(z)\Gamma(dz)ds \notag\\
  & = \int_0^T\int_\Rd |\tf g(s)(\xi)|^2 \mu_s^Z(d\xi)ds. \label{eq:normP0Z} 
\end{align}

%Assumptions on \Lambda
Following \cite{conusdalang}, we describe the assumptions on deterministic functions that may be integrated with respect to the martingale measure $M^Z$. These are as follows.

\begin{description}
\item{{\bf (A3)}}
$t\mapsto\Lambda(t)$ is a deterministic function with values in $\mathcal{S}'_r(\Rd)$; the mapping $(t,\xi)\mapsto\tf\Lambda(t)(\xi)$ is measurable and
  \[ \int_0^T \sup_{\eta\in\Rd} \int_{\Rd}  |\tf\Lambda(s)(\xi+\eta)|^2 \mu(d\xi) ds< \infty. \]
 \item{{\bf(A4)}} Let $\phi$ denote a nonnegative function in $\mathcal{C}^\infty_0(\Rd)$, with support included in the unit ball of $\Rd$, satisfying $\int_\Rd \phi(x)dx=1$. For all such $\phi$ and all $0\leq a\leq b\leq T$, we have
  \[ \int_a^b (\Lambda(s)\ast\phi)(x) ds \in \mathcal{S}(\Rd) \]
  and
  \[ \int_{\Rd}\int_a^b |(\Lambda(s)\ast\phi)(x)| ds dx < \infty. \]
\item{{\bf (A5)}} $t\mapsto\tf\Lambda(t)$ is as in {\bf (A3)} and
  \[ \lim_{h\downarrow0} \int_0^T \sup_{\eta\in\Rd}\int_{\Rd}\sup_{s<r<s+h}|\tf\Lambda(r)(\xi+\eta)-\tf\Lambda(s)(\xi+\eta)|^2\mu(d\xi)\ ds = 0. \]
\end{description}  
Notice that {\bf (A3)} implies 
  \begin{equation}
		\int_0^T \int_{\Rd} |\tf\Lambda(s)(\xi)|^2 \mu(d\xi)ds < \infty.
		\label{eq:classiccondition}
	\end{equation} 

In \cite[Theorem 3.1]{conusdalang} it is proved that under assumptions {\bf (A1)}, {\bf (A2)}, {\bf (A3)}, and either {\bf (A4)} or {\bf (A5)}, 
$\Lambda\in\mathcal{P}_{0,Z}$ and that the stochastic integral $((\Lambda\cdot M^Z)_t;\; t\in[0,T])$ is well-defined as a real-valued square-integrable martingale. Moreover,  
\begin{equation}
	\E\big[(\Lambda\cdot M^Z)_t^2\big] = \int_0^t\int_\Rd |\tf\Lambda(s)(\xi)|^2\mu_s^Z(d\xi)ds = \|\Lambda\|_{0,Z}^2.
	\label{eq:normstoint}
\end{equation} 
Using Assumption {\bf (A2)} and the identities $\nu^Z_s(\Rd) = \gamma^Z_s(0) = \E[Z(s,0)^2]$, one can obtain the following upper bound:  

\begin{align}
  \E\big[(\Lambda\cdot M^Z)^2_t\big] 
  & = \int_0^t\int_{\Rd}|\tf\Lambda(s)(\xi)|^2\mu_s^Z(d\xi) ds \notag\\
  %& = \int_0^t\int_{\Rd}\int_{\Rd} |\tf\Lambda(t-s)(\xi+\eta)|^2\mu(d\xi)\nu^Z_s(d\eta)do s \notag\\
  & \leq \int_0^t\nu^Z_s(\Rd)\sup_{\eta\in\Rd}\int_{\Rd} |\tf\Lambda(s)(\xi+\eta)|^2\mu(d\xi) ds \notag\\
  & = \int_0^t\E[Z(s,0)^2]\sup_{\eta\in\Rd}\int_{\Rd} |\tf\Lambda(s)(\xi+\eta)|^2\mu(d\xi) ds, \label{eq:boundCDI}
\end{align}
Owing to {\bf (A1)}, this yields 
\begin{equation}
\label{1.1}
 \E\big[(\Lambda\cdot M^Z)^2_t\big]\le C
	\int_0^t\sup_{\eta\in\Rd}\int_{\Rd} |\tf\Lambda(s)(\xi+\eta)|^2\mu(d\xi)ds.
\end{equation}

%Pathwise integral
Let $t\mapsto\psi(t)$ be a deterministic functions with values in $\mathcal{S}'_r(\Rd)$. Assume that $\psi\in L^2([0,T];L^1(\Rd))$, that is, satisfying $\int_0^T\left(\int_{\Rd} \vert \psi(s,z)\vert dz \right)^2 ds < \infty$. For a stochastic process $Z$ satisfying the conditions {\bf (A1)}, {\bf (A2)}, and following again \cite{conusdalang}, we introduce the norm
 \[ \|\psi\|_{1,Z}^2 := \E\bigg[\int_0^T\bigg(\int_\Rd \psi(s,z)Z(s,z) dz\bigg)^2ds\bigg]. \]
Proceding as in the derivation of \eqref{eq:normP0Z}, we obtain
\begin{align}
\label{1.2}
  \|\psi\|_{1,Z}^2
  & := \E\bigg[\int_0^T\bigg(\int_\Rd \psi(s,z)Z(s,z) dz\bigg)\bigg(\int_\Rd \psi(s,y)Z(s,y)dy \bigg) ds\bigg] \notag\\
  & = \int_0^T  \int_\Rd  \int_\Rd  \psi(s,z)g(s,y)\gamma_s^Z(z-y)\ dy dz ds\notag\\
  & = \int_0^T  \int_\Rd |\tf \psi(s)(\eta)|^2 \nu_s^Z(d\eta) ds.
\end{align}  
The closure of the space $\mathcal{E}$ with respect to the norm $\Vert\cdot\Vert_{1,Z}$ is denoted by $\mathcal{P}_{1,Z}$. 

In order to give a rigorous meaning to pathwise convolutions, some additional assumptions are needed. These are the following.

\begin{description}
\item{{\bf(A6)}}
The mapping $t\mapsto\Lambda(t)$ is a deterministic function with values in $\mathcal{S}_r^\prime(\Rd)$ and satisfies 
 \begin{equation*}\int_0^T \sup_{\eta\in\Rd} |\tf\Lambda(s)(\eta)|^2 ds < \infty. 
 \end{equation*}
 \item{{\bf(A7)}}
The mapping $t\mapsto\tf\Lambda(t)$ is a deterministic function with values in $\mathcal{S}_r^\prime(\Rd)$ and such that
    \[ \lim_{h\downarrow0} \int_0^T \sup_{\eta\in\Rd}\sup_{s<r<s+h}|\tf\Lambda(r)(\eta)-\tf\Lambda(s)(\eta)|^2  ds = 0. \]
\end{description}

Note that these two conditions coincide respectively with {\bf(A3)} and {\bf(A5)} if $\mu=\delta_0$. 

Assume {\bf(A1)}, {\bf(A2)}, {\bf(A6)} and either {\bf(A4)} or {\bf(A7)}.
In \cite[Proposition 3.4]{conusdalang} it is proved that
\begin{equation}
\label{1200}
\int_0^t\int_\Rd \Lambda(s,z)Z(s,z) dzds, \ t\in[0,T], 
\end{equation}
defines a stochastic process with  values in $L^2(\Omega)$. 
In addition, from \eqref{1.2} and {\bf(A1)}, {\bf(A2)}, it follows that 
\begin{align}
\label{1.3}
  \Vert \Lambda\Vert_{1,Z}^2=& \E\Bigg[\bigg(\int_0^t\int_{\Rd} \Lambda(s,x)Z(s,x) dx ds\bigg)^2\Bigg] \nonumber\\%\int_0^t\int_\Rd |\tf\Lambda(s)(\eta)|^2\nu^Z_s(d\eta)do s \notag\\
        &\leq \int_0^t \nu^Z_s(\Rd) \sup_{\eta\in\Rd} |\tf\Lambda(s)(\eta)|^2 ds \notag\\
       % & = \int_0^t \E[Z(s,0)^2] \sup_{\eta\in\Rd} |\tf\Lambda(s)(\eta)|^2do s, 
       & \leq C  \int_0^t \sup_{\eta\in\Rd} |\tf\Lambda(s)(\eta)|^2 ds.
\end{align}

Frequently, we will use the notation
\begin{equation}
  J_1(s) := \sup_{\eta\in\Rd}\int_{\Rd} |\tf\Lambda(s)(\xi+\eta)|^2\mu(d\xi)
\label{eq:J1}
\end{equation}
and
\begin{equation}
  J_2(s) := \sup_{\eta\in\Rd} |\tf\Lambda(s)(\eta)|^2,
\label{eq:J2}
\end{equation}
$s\in[0,T]$.

If the assumptions {\bf(A3)}, {\bf(A6)}, respectively, are satisfied, then 
\begin{equation}
	\int_0^T J_1(s) ds < \infty, \ \int_0^T J_2(s) ds < \infty,  
\label{eq:2.17'''}
\end{equation}
respectively.

%%%%%last notations and facts

Throughout the article, we will refer extensively to the Hilbert space $L^2(\Omega;{\mathcal{H}_T})=L^2([0,T]\times\Omega;\mathcal{H})$. In \cite[Proposition 2.6]{dalangquer}, it is proved that 
$\mathcal{P}_0\subseteq L^2(\Omega;{\mathcal{H}_T})$. Then if $g\in \mathcal{P}_0$, $\|g\|_0=\|g\|_{L^2(\Omega;{\mathcal{H}_T})}$ with $\Vert\cdot\Vert_0$ defined as in \eqref{eq:normP01}.

In this article, we will use the theory of Malliavin calculus based on the Gaussian process $F=(F(\phi);\; \phi\in{\mathcal{H}_T})$ (see \cite{nualart}). For this we need to guarantee that $F$ is an isonormal Gaussian process and also to describe its associated abstract Wiener space. 

By Lemma \ref{lem:P0=HT}, the expression \eqref{eq:FandM} holds for any $\phi\in {\mathcal{H}_T}$. This yields that $F$ is an isonormal process (see \cite[Definition 1.1.1]{nualart}). For the description of the abstract Wiener space, it is useful to identify the stochastic process $F$ with a $\hac$--valued cylindrical Wiener process, as follows.
As it is shown in \cite{dalangfrangos}, by an approximation procedure we define  $W_t(\phi)=F(1_{[0,t]}\phi)$, $t\in[0,T]$, $\phi\in\hac$. Consider a complete orthonormal system
(CONS) of $\hac$ that we denote by $(e_k)_{k\in\N}$. Then, 
\begin{equation*}
W=\{W^k(t):=W_t(e_k), t\in[0,T], k\in\N\}
\end{equation*}
defines a sequence of independent standard Brownian motions.
Conversely, the process $(F(\phi)= \sum_{k\in\N}Ê\int_0^T \langle\phi(t),e_k\rangle_{\hac} dW^k(t), \phi\in\hac_T)$ is an isonormal Gaussian process.

Let $(\bar\Omega,\bar{\mathcal{G}}, \bar \mu)$ be the canonical space of a standard real-valued Brownian motion on $[0,T]$.
With the equivalence shown before, we can identify the canonical probability space of $F$ with that of a sequence of independent standard Brownian motions $(\Omega,\mathcal{G},\mathbb{P})=(\bar\Omega^{\N}, \bar{\mathcal{G}}^{\otimes\N}, \bar{\mu}^{\otimes\N})$. This will be the
underlying probability space in this work. 

Consider the Hilbert space $\mathbb{H}$ consisting of sequences $(h^k)_{k\in\N}$ of functions $h_k:[0,T]\rightarrow \R$
which are absolutely continuous with respect to the Lebesgue measure and such that $\sum_{k\in\N}\int_0^T|\dot h^k(s)|^2 ds<\infty$, where $\dot h^k$ refers to the derivative of $h^k$ defined almost everywhere. There is an isometry between the spaces $\mathbb{H}$ and $\hac_T$, as follows. Let $h\in\mathbb{H}$. Then $h=\sum_{k\in\N}h^ke_k$, where $h^k(t)=\int_0^t \dot{h}^k(s) ds$ for all $k\in\N$. For any $t\in[0,T]$, set 
$\bar{h}(t)=\sum_{k\in\N} \dot{h}^k(t)e_k$. Clearly, $\bar{h}\in\hac_T$ and $\|h\|_{\mathbb{H}} = \|\bar{h}\|_{\hac_T}$. The triple $(\Omega,\mathbb{H},\mathbb{P})$ is the abstract Wiener space that we shall use as framework for the Malliavin calculus.
\smallskip

%%%%%%%%%%
%%%%%%%%%% END SECTION PRELIMINARIES

%%%%%%%%%%
%%%%%%%%%% SECTION H-STOCHASTIC INTEGRALS
%%%%%%%%%%%%%%

\section{Hilbert space valued stochastic integrals}
\label{s2}
In this section, we develop an extension of the integrals introduced in \cite{conusdalang} to Hilbert space valued integrands. For similar results in the setting of \cite{dalang}, we refer the reader to \cite{quersanz}.

Let $\caA$ be a separable real Hilbert space with inner-product and norm denoted by $\langle\cdot,\cdot\rangle_\caA$ and $\|\cdot\|_\caA$ respectively. In the sequel, $(a_k)_{k\in\N}$ 
will denote  a complete orthonormal system of $\caA$. Let $\{Z(t,x), t\in[0,T], x\in\Rd\}$ be a $\caA$-valued stochastic process satisfying the following properties similar to {{\bf(A1)}}, {{\bf(A2)}}:
\medskip

\noindent{\bf(A8)}  
$Z$ is a $\caA$-valued predictable stochastic process satisfying
$$\sup_\tx \E[\|Z(t,x)\|_\caA^2] < \infty.$$ 

\noindent{\bf(A9)} $Z$ has a {\it spatial stationary covariance function} coordinatewise. That is, for all $k\in\N$ and $(t,x)\in [0,T]\times \Rd$,
    \[ \E[Z_k(t,x)Z_k(t,x+y)] = \E[Z_k(t,0)Z_k(t,y)] =: \gamma_{k,t}^Z(y), \]
    where $Z_k(t,x):= \langle Z(t,x),a_k\rangle_{\caA}$.

For each $k\in\N$, the stochastic process
$$
M^{Z_k}_t(A) := \int_0^t\int_A Z_k(s,z)M(ds,dz), t\in[0,T], \ A\in\mathcal{B}_b(\Rd),
$$
defines a real-valued worthy martingale measure. 

According to \cite{quersanz}, %\cite[p.\ 5]{quersanz}, 
$$
 M^Z_t(A) := \sum_{k\in\N} M^{Z_k}_t(A)a_k, \ t\in[0,T], \ A\in\mathcal{B}_b(\Rd),
 $$
defines a $\caA$--valued worthy martingale measure and by construction,\break $\langle M_t^Z(A), a_k\rangle_\caA = M_t^{Z_k}(A)$. Moreover, 
\begin{equation*}
  \E\big[\|M_t^Z(A)\|_\caA^2\big] \leq C_{t,A}\sup_\tx \E\big[\|Z(t,x)\|_\caA^2\big],
\end{equation*}
(see \cite[p.\ 5]{quersanz}).
%Here and in the sequel, the symbols ``$\cdot$", ``$\ast$" refer to the temporal and the spatial arguments, respectively. 

For a $\caA$--valued stochastic process $g$, we set $g_k=\langle g,a_k\rangle_{\caA}$, and define
$$
\Vert g\Vert_{0,\caA}^2 := \sum_{k\in\N}\Vert g_k\Vert_0^2.
$$
Also, for a deterministic function $\phi$, we define
\begin{equation*}
\|\phi\|_{0,Z,\caA}^2 := \sum_{k\in\N} \|\phi\|_{0,Z_k}^2 ,\quad  \|\phi\|_{+,Z,\caA}^2 := \sum_{k\in\N} \|\phi\|_{+,Z_k}^2,\quad  \|\phi\|_{1,Z,\caA}^2 := \sum_{k\in\N} \|\phi\|_{1,Z_k}^2.
\end{equation*}
In each one of these definitions we are implicitely assuming that the right-hand side of every expression is well-defined. Let $\mathcal{E}_\caA$ be the set of $\caA$-valued simple stochastic processes. That is, processes with a similar expression as in \eqref{eq:simpleprocess}, where $X_j$, $j=1,\ldots, m$, are $\caA$-valued  random variables. Then, we denote by $\mathcal{P}_{0,Z,\caA}$, $\mathcal{P}_{+,Z,\caA}$ and $\mathcal{P}_{1,Z,\caA}$ the completions of $\mathcal{E}$ with respect to $\|\cdot\|_{0,Z,\caA}$, $\|\cdot\|_{+,Z,\caA}$  and $\|\cdot\|_{1,Z,\caA}$, respectively. 

Note that all the norms $\|\cdot\|_{0,Z,\caA}, \|\cdot\|_{+,Z,\caA}$ and $\|\cdot\|_{1,Z,\caA}$ do not depend on the choice of the CONS although assumption ${\bf (A9)}$ might do so. Indeed, for the dense subset of $\mathcal{P}_{0,Z,\caA}$ for which the norm $\|\cdot\|_{0,Z,\caA}$ can be written as in \eqref{eq:normP0Z1}, one can easily verify that $\sum_{k\in\N} \|\cdot\|_{0,\langle Z,a_k\rangle_\caA}^2 = \sum_{k\in\N} \|\cdot\|_{0,\langle Z,a'_k\rangle_\caA}^2$, where $(a_k)_{k\in\N}$ and $(a'_k)_{k\in\N}$ are two CONS of $\caA$. By density, this equality then holds for all elements in $\mathcal{P}_{0,Z,\caA}$.

%%%%%
%%%%% DOS TEOREMES D'EXTENSIO 
%%%%%

\begin{theorem}
\label{t2.1}
Let $\{Z(t,x), (t,x)\in[0,T]\times \Rd\}$ be a $\caA$--valued stochastic process satisfying conditions {\bf(A8)}, {\bf(A9)}. Let  $t\mapsto \Lambda(t)$ be a deterministic function taking values in the space $\mathcal{S}^{\prime}_r(\Rd)$. We suppose that  {\bf(A3)} and either {\bf(A4)} or {\bf(A5)} are satisfied. Then 
$\Lambda\in \mathcal{P}_{0,Z,\caA}$ and the stochastic integral $\{(\Lambda\cdot M^Z)_t, t\in[0,T]\}$ is well-defined as a $\caA$-valued square integrable process. Moreover,
\begin{align}
\label{2.1}
  \E\big[\|(\Lambda\cdot M^Z)_t\|_\caA^2\big]&= \int_0^t\int_\Rd |\tf\Lambda(s)(\xi)|^2 \mu_s^Z(d\xi)ds = \|\Lambda\|_{0,Z,\caA}^2\nonumber\\
& \le  \int_0^t \E[\|Z(s,0)\|_\caA^2] \sup_{\eta\in\Rd}\int_\Rd |\tf\Lambda(s)(\xi+\eta)|^2 \mu(d\xi)ds\nonumber\\
&\le C   \int_0^t \sup_{\eta\in\Rd}\int_\Rd |\tf\Lambda(s)(\xi+\eta)|^2 \mu(d\xi)ds
 \end{align}
\end{theorem}

\noindent{\it Proof.} From \cite[Theorem 3.1]{conusdalang} we know that $\Lambda\in\mathcal{P}_{0,Z_k}$ and also that 
$\{(\Lambda\cdot M^{Z_k})_t, t\in[0,T]\}$ is well-defined for any $k\in\N$.
In addition,
\begin{equation*}
 \E\big[\|(\Lambda\cdot M^{Z_k})_t\|^2\big]=\Vert \Lambda\Vert_{0,Z_k}^2.
 \end{equation*}
 
 The proof of $\Lambda\in \mathcal{P}_{0,Z,\caA}$ follows the same arguments as in \cite[Theorem 3.1]{conusdalang}.
 Firstly, we check that $\Lambda_n$ (defined similar to \eqref{eq:Lambda_n} by $\Lambda_n(t):=\Lambda(t)\ast\zeta_n$) belongs to $\mathcal{P}_{0,Z,\caA}$, and then that
 \begin{equation*}
 \lim_{n\to\infty}\Vert\Lambda_n-\Lambda\Vert_{0,Z,\caA}^2=0.
 \end{equation*}
 The arguments of \cite[Theorem 3.1]{conusdalang} can be adapted by using {\bf{(A3)}}, {\bf{(A8)}} and the following remark:
 For any $\phi:[0,T]\to\mathcal{S}(\Rd)$,
 \begin{equation*}
 \Vert\phi\Vert^2_{0,Z_k}\le \int_0^TE\big[|Z^k(s,0)|^2\big] \sup_{\eta\in\Rd}\int_{\Rd} \mu(d\xi)|\tf \phi(s)(\xi+\eta)|^2.
 \end{equation*}
 It follows that
 \begin{align*}
 \Vert\phi\Vert^2_{0,Z,\mathcal{A}}&=\sum_{k\in\N} \Vert\phi\Vert^2_{0,Z_k}\\
 &\le \int_0^TE\big[\Vert Z(s,0)\Vert_{\mathcal{A}}^2\big] \sup_{\eta\in\Rd}\int_{\Rd} \mu(d\xi)|\tf \phi(s)(\xi+\eta)|^2.
 \end{align*}
 
 For any $t\in[0,T]$, define
 \begin{equation*}
 (\Lambda\cdot M^Z)_t=\sum_{k\in\N} (\Lambda\cdot M^{Z_k})_t a_k.
 \end{equation*}
 Clearly,
 \begin{equation*}
  \E\big[\|(\Lambda\cdot M^Z)_t\|_\caA^2\big]=\sum_{k\in\N} \E\big[\|(\Lambda\cdot M^{Z_k})_t\|^2\big]= \sum_{k\in\N}\Vert\Lambda\Vert_{0,Z_k}^2
  = \|\Lambda\|_{0,Z,\caA}^2.
  \end{equation*}
  The estimates \eqref{2.1} follows from \eqref{eq:boundCDI} applied to each stochastic integral $\Lambda\cdot M^{Z_k}$, $k\in\N$, along with {\bf{(A8)}}.
  
  \cqd
  
  By using similar arguments as in the proof of  \cite[Proposition 3.4]{conusdalang}, one can also give an extension of the pathwise integral to $\caA$--valued stochastic processes. For this, it is worth noticing that
  for any $\phi:[0,T]\to\mathcal{S}(\Rd)$,
    \begin{align*}
    \Vert\phi\Vert_{1,Z,\mathcal{A}}^2&=\sum_{k\in\N}\Vert\phi\Vert_{1,Z^k}\\
    &\le \int_0^t \E[\Vert Z(s,0)\Vert^2_{\mathcal{A}}]\sup_{\eta\in\Rd}|\tf\phi(s)(\eta)|^2 ds.
    \end{align*}
    The extension reads as follows.
  
    \begin{theorem}
 \label{t2.2}
 Let $\{Z(t,x), (t,x)\in[0,T]\times \Rd\}$ be a stochastic process as in Theorem \ref{t2.1}. Let  $t\mapsto \Lambda(t)$ be a deterministic function taking values in the space $\mathcal{S}^{\prime}_r(\Rd)$. We suppose that  {\bf(A6)} and either {\bf(A4)} or {\bf(A7)} are satisfied. Then %for any $x\in\Rd$,
 \begin{align*}
 &\int_0^t\int_{\Rd}\Lambda(s,z) Z(s,z) dzds
:=\sum_{k\in\N}\left( \int_0^t\int_{\Rd}\Lambda(s,z) Z_k(s,z) dzds\right)a_k,
 \end{align*}
 $t\in[0,T]$, defines a stochastic process with values in $L^2(\Omega;\caA)$. Moreover,

 \begin{align}
 \label{2.3}
  \E\left[\left\|\int_0^t\int_\Rd \Lambda(s,z)Z(s,z) dzds\right\|^2_\caA\right]
  & = \sum_{k\in\N} \E\left[\left(\int_0^t\int_\Rd \Lambda(s,z)Z_k(s,z) dzds\right)^2\right]\nonumber\\
  &=\Vert\Lambda\Vert_{1,Z,\caA}^2 \notag\\
  & \leq \sum_{k\in\N} \int_0^t \E[Z_k(s,0)^2] \sup_{\eta\in\Rd} |\tf\Lambda(s)(\eta)|^2 ds \notag\\
  %& = \int_0^t \E[\|Z(s,0)\|_\caA^2] \sup_{\eta\in\Rd} |\tf\Lambda(s)(\eta)|^2 do s\notag\\
%  & \leq \int_0^t \sup_{y\in\Rd}\E[\|Z(s,y)\|_\caA^2] \sup_{\eta\in\Rd} |\tf\Lambda(s)(\eta)|^2 do s\\
 % & \leq \int_0^t \sup_{(r,y)\in[0,s]\times\Rd}\E[\|Z(r,y)\|_\caA^2] \sup_{\eta\in\Rd} |\tf\Lambda(s)(\eta)|^2 do s \notag\\
  & \leq C \int_0^t \sup_{\eta\in\Rd} |\tf\Lambda(s)(\eta)|^2 ds.   
\end{align}
\end{theorem}

 %%%%%%%%%%
%%%%%%%%%% END SECTION H-STOCHASTIC INTEGRALS

%%%%%%%
%%%%%%%EQUIVALENCE OF STOCHASTIC INTEGRALS

\section{Equivalence of stochastic integrals}
\label{s3}

In this section, we consider a particular case of integrands described as follows. Let $Z$ be a stochastic process satisfying {\bf(A1)}, {\bf(A2)}. Let $\Lambda:[0,T]\to\mathcal{S}_r^\prime(\Rd)$. We are interested in stochastic processes which are obtained as the limit in the topology of $\mathcal{P}_0$ of a sequence
\begin{equation*}
\Phi_{t,x}^n:=\Lambda_n(t-\cdot,x-\ast)Z(\cdot,\ast), n\in\N,
\end{equation*}
where $(t,x)\in[0,T]\times \R^d$ is fixed and $\Lambda_n(t):=\Lambda(t)\ast\zeta_n$ as in \eqref{eq:Lambda_n}.

For this class of integrands, later denoted by $\Phi_{t,x}$ (and also at some places by $\Lambda(t-\cdot,x-\ast)Z(\cdot,\ast)$, by an abuse of language), we prove that
the integrals in the Conus-Dalang sense (\cite{conusdalang}) and with respect to the $\hac$-valued cylindrical Wiener process $(W_t, t\in[0,T])$ (see for instance \cite{dapratozabczyk}), coincide with the divergence operator (also termed Skorohod integral) of Malliavin Calculus (see \cite[Section 1.3]{nualart}). 
For this, we need further insight on the relationships between the spaces $\mathcal{P}_0$, $\mathcal{P}_{0,Z}$ and $L^2(\Omega,{\mathcal{H}_T})$ introduced in Section \ref{s1}. 

%Let throughout this section and the remaining paper $\Lambda_n$ denote the regularization of $\Lambda$ in \eref{eq:Lambda_n}. Moreover, for the sake of brewity, let for any fixed $\tx$, $\Lambda_nZ := \Lambda_n(t-\cdot,x-\ast)Z(\cdot,\ast)$ and $\Lambda_nDZ := \Lambda_n(t-\cdot,x-\ast)DZ(\cdot,\ast)$ where $DZ$ denotes the Malliavin Derivative of $Z(\cdot,\ast)$.\\
%At first we note that in \cite[Proposition 2.6(b)]{dalangquer} it was shown that an element $g\in\mathcal{P}_0$ also belongs to $L^2(\Omega,{\mathcal{H}_T})$. 

We notice that, for deterministic elements $\phi\in\mathcal{P}_0$, which are Schwartz functions in the spatial argument and a process $Z$ satisfying {\bf(A1)}, {\bf(A2)}, 
\begin{equation}
	\|\phi Z\|_0 = \|\phi\|_{0,Z} = \|\phi Z\|_{L^2(\Omega;{\mathcal{H}_T})}.
\label{3.1}
\end{equation} 
%In fact, these equalities hold more generally when $g$ is a function in $\mathcal{P}_{0,Z}$. Note that we assume $Z$ to be a random field (which permits a pointwise interpretation as a random variable in $L^2(\Omega)$ for every given $\tx$). \\
%Also, the norms for $\caA$-valued processes introduced in Section \ref{sec:hilbert} satisfy such isometries. If $g$ is a Schwartz function and $Z$ a process satisfying \ref{itm:DZ^2} and \ref{itm:DZstat} then the previous equality holds, i.e.
%\begin{equation}
 % \|g Z\|_{0,\caA} = \|g\|_{0,Z,\caA}. 
%\label{eq:equalofnormsHS}
%\end{equation}

%The following lemma gives establishes some more relationships between elements of $\mathcal{P}_0$ and those of $\mathcal{P}_{0,Z}$ where the first two assertions were already shown in \cite[Theorem 3.1]{conusdalang}. 
To simplify the notation, for any $(t,x)\in[0,T]\times\Rd$, we write $[\Lambda_nZ]^{t,x}$ to denote the stochastic process $(\Lambda_n(t-\cdot,x-\ast)Z(\cdot,\ast))$, where ``$\cdot$" and ``$\ast$" denote the time and space arguments, respectively.%, where $\Lambda_n$ is defined in \eqref{eq:Lambda_n}
\begin{lemma}
\label{l3.1}
Let $Z$ be a stochastic process satisfying the hypotheses {\bfÊ(A1)}, {\bfÊ(A2)}. 
 Let $t\in[0,T]\mapsto \Lambda(t)$ be a function satisfying the assumptions {\bfÊ(A3)} and either {\bfÊ(A4)} or {\bfÊ(A5)}. 
Fix $(t,x)\in[0,T]\times\Rd$. Then, for all $n\in\N$,
\begin{enumerate}
	\item $\Lambda_n, \Lambda\in\mathcal{P}_{0,Z}$, 
	\item $[\Lambda_nZ]^{t,x}\in\mathcal{P}_0$,
	\item\label{itm:Phi} The sequence $([\Lambda_nZ]^{t,x})_{n\in\N}$ converges in $\mathcal{P}_0$ to an element $\Phi_{t,x}\in\mathcal{P}_0$, and  
	\begin{equation}
		\|\Phi_{t,x}\|_0 = \|\Lambda(t-\cdot,x-\ast)\|_{0,Z}.
	\label{3.2}
	\end{equation}
	Moreover, $\Phi_{t,x}\in L^2(\Omega; \hac_T)$ and 
	\begin{equation*}
		\|\Phi_{t,x}\|_0 = \|\Lambda(t-\cdot,x-\ast)\|_{L^2(\Omega;\hac_T)}.
	%\label{3.2}
	\end{equation*}
\end{enumerate}	
\end{lemma}	
\noindent{\it Proof.} The assertions of part (1) are shown in \cite[Theorem 3.1]{conusdalang}.
%proof of (ii)
The second part in shown by a similar method as in \cite[Theorem 3.1]{conusdalang}. In fact, consider either approximation of $\Lambda_n$ by simple functions $(\Lambda_{n,m})_{m\in\N}$ given in the proof of this theorem. Then one shows using the definition of the norm $\Vert\cdot\Vert_+$
\begin{align}
\label{3.3}
  & \|\Lambda_{n,m}(t-\cdot,x-\ast)Z(\cdot,\ast)\|^2_+ \notag\\
  & = \E\Big[\int_0^t\int_\Rd\int_\Rd |\Lambda_{n,m}(t-s,x-z)| |\Lambda_{n,m}(t-s,x-y+z)| \nonumber\\
  &\quad \times |Z(s,y)||Z(s,y-z)| dy\Gamma(dz)ds\Big] \notag\\
  & \leq \sup_{(r,y)\in[0,T]\times\Rd}\E[Z(r,y)^2] \nonumber\\
  &\quad \times \int_0^t\int_\Rd \big(|\Lambda_{n,m}(t-s,x-\cdot)|\ast|\tilde{\Lambda}_{n,m}(t-s,x-\cdot)|\big)(z) \Gamma(dz)ds.
\end{align}
Since $\Lambda_{n,m}(t)\in\mathcal{S}(\Rd)$ and $\Gamma$ is a tempered measure, we can use Leibniz' rule (\cite[Exercise 26.4]{treves}) to shown that the last expression is finite. This shows that $[\Lambda_{n,m}Z]^{t,x}\in\mathcal{P}_+$. Then we evaluate the difference $\|(\Lambda_{n,m}(t-\cdot,x-\ast)-\Lambda_n(t-\cdot,x-\ast))Z(\cdot,\ast)\|_0^2$ in the same ways (depending on whether we suppose {\bfÊ(A4)} or {\bfÊ(A5)}) as in the proof of \cite[Theorem 3.1]{conusdalang}, and show that it goes to zero. This proves part 2.
%proof of (iii)

For $n, m\in\N$, we have
%In order to prove (iii), we show that $(\Lambda_nZ)\nN$ is a Cauchy sequence in $\mathcal{P}_0$. For this we show that $(\Lambda_n)\nN$ is a $\|\cdot\|_{0,Z}$-Cauchy sequence:
\begin{align}
\label{3.4}
  \|\Lambda_n & (t-\cdot,x-\ast)-\Lambda_m(t-\cdot,x-\ast)\|_{0,Z}^2 \nonumber\\
              & = \int_0^t\int_\Rd |\tf\Lambda_n(t-s)(\xi) - \tf\Lambda_m(t-s)(\xi)|^2 \mu^Z_s(d\xi)ds\nonumber\\
              & = \int_0^t\int_\Rd |\tf\Lambda(t-s)(\xi)|^2 |\tf\zeta_n(\xi)-\tf\zeta_m(\xi)|^2 \mu^Z_s(d\xi)ds.
\end{align}
By bounded convergence, this converges to zero as $n,m\to\infty$. Using \eqref{3.1}, we conclude that 
$([\Lambda_nZ]^{t,x})_{n\in\N}$ is a Cauchy sequence in $\mathcal{P}_0$. Let us denote by $\Phi_{t,x}$ its limit. 
Using similar computations as in \eqref{3.4} with $\Lambda_m$ replaced by $\Lambda$, and since $\Lambda\in\mathcal{P}_{0,Z}$,
we have,
\begin{equation*}
\lim_{n\to\infty} \|\Lambda_n(t-\cdot,x-\ast)\|_{0,Z} 
                            = \|\Lambda(t-\cdot,x-\ast)\|_{0,Z}.
\end{equation*}
Thus,
\begin{equation*}
  \|\Phi_{t,x}(\cdot,\ast)\|_0 %= \lim_{n\to\infty} \|\Lambda_n(t-\cdot,x-\ast)Z(\cdot,\ast)\|_0
                               = \lim_{n\to\infty} \|\Lambda_n(t-\cdot,x-\ast)\|_{0,Z} \notag\\
                               = \|\Lambda(t-\cdot,x-\ast)\|_{0,Z}. 
\end{equation*}
Since $\mathcal{P}_0\subseteq L^2(\Omega;\hac_T)$ (see \cite[Proposition 2.6]{dalangquer}), we conclude that $\Phi_{t,x}\in L^2(\Omega;\hac_T)$. \phantom{mmmmm} \cqd
%%%%%%new version ends here (Corollary 4.2 from AndrŽ)
%%%%%%

The preceding lemma admits easily an extension to Hilbert space valued stochastic processes. Next, we consider a particular example of such an extension for processes that are Malliavin derivatives.

\begin{lemma}
\label{l3.2}
The function $\Lambda$ and the stochastic process $Z$ are as in Lemma \ref{l3.1}. Moreover,  we assume that that $Z(t,x)\in\D^{1,2}$ for all $(t,x)\in[0,T]\times\Rd$ and that the $\hac_T$--valued stochastic process 
%$(DZ(t,x),(t,x)\in[0,T]\times \R^d)$
$DZ$ satisfies {\bf(A8)}, {\bf(A9)}. Fix $(t,x)\in[0,T]\times \R^d$, Then, by setting 
\begin{equation*}
[\Lambda_nDZ]^{t,x}:=\left(\Lambda_n(t-s,x-y)DZ(s,y), (s,y)\in[0,T]\times \R^d\right),
\end{equation*}
the following holds.
\begin{enumerate}
\item $\Lambda_n, \Lambda\in\mathcal{P}_{0,DZ,{\mathcal{H}_T}}$,
\item $[\Lambda_nDZ]^{t,x}\in\mathcal{P}_{0,{\mathcal{H}_T}}$ 
\item The sequence $([\Lambda_nDZ]^{t,x})_{n\in N}$ converges in $\mathcal{P}_{0,{\mathcal{H}_T}}$ to a $\hac_T$--valued stochastic process 
$\Phi^{(1)}_{t,x}:=\{\Phi_{t,x}^{(1)}(s,y), (s,y)\in[0,T]\times \R^d\}$ such that for any $(t,x)\in[0,T]\times \R^d$,
	\begin{equation} 
		\|\Phi_{t,x}^{(1)}\|_{0,{\mathcal{H}_T}} = \|\Lambda(t-\cdot,x-\ast)\|_{0,DZ,\mathcal{H}_T}.
		\label{3.40}
	\end{equation}	
	\item  $\Phi^{(1)}_{t,x}=D \Phi_{t,x}$, where $\Phi_{t,x}$ is the process defined in part 3 of Lemma \ref{l3.1}.
	%The Malliavin derivative of the limit $\Phi_{t,x}$ obtained in Lemma \ref{lem:LambdaZ} is almost surely equal to $\Phi_{t,x}^{(1)}$.
\end{enumerate}
\end{lemma}

{\it Proof.}
Statement 1 can be shown as in \cite[Theorem 3.1]{conusdalang} with the tools provided in Section \ref{s2}. 
For the proof of part 2, we follow similar computations as in \eqref{3.3} to obtain
\begin{align*}
&\Vert \Lambda_{m,n}(t-\cdot,x-\ast)DZ_{\cdot,\ast}(t,x)\Vert_{+,\hac_T}^2 \le C\sup_{(r,y)\in[0,T]\times\Rd} \E\big[\|DZ(r,y)\|^2_{\mathcal{H}_T}\big]\\ 
&\qquad\times\int_0^t\int_\Rd \big(|\Lambda_{n,m}(t-s,x-\cdot)|\ast|\tilde{\Lambda}_{n,m}(t-s,x-\cdot)|\big)(z) \Gamma(dz)ds.
\end{align*}
By the same arguments as in Lemma \ref{l3.1}, this last expression is finite.

Similarly as for the statement 3 of Lemma \ref{l3.1}, we prove that $([\Lambda_nDZ]^{t,x}, n\in \N)$ is a Cauchy sequence in the norm $\Vert\cdot\Vert_{0,DZ,\mathcal{H}_T}$ and that its limit $\Phi^{(1)}$ satisfies \eqref{3.40}. 

As for part 4, we notice that by Lemma \ref{l3.1}, the sequence $(\Lambda_n(t-\cdot,x-\ast)Z(\cdot,\ast), n\in\N)$ converges in $L^2(\Omega; \hac_T)$ to a random vector 
$\Phi_{t,x}$.  Moreover, $D(\Lambda_n(t-\cdot,x-\ast)Z)= \Lambda_n(t-\cdot,x-\ast)DZ$. Hence by part 3, the sequence $\left(D(\Lambda_n(t-\cdot,x-\ast)Z)\right)_{n\in\N}$  converges in $L^2(\Omega; \hac_T^{\otimes 2})$. Since $D$ is a closed operator, we conclude using again part 3. 

\cqd

Let $\Lambda$ and $Z$ be as in Lemma \ref{l3.1} and $(\bar{e}_k)_{k\in\N}$ be a CONS of $\hac_T$. For any $k\in\N$, the real-valued stochastic process $D^{\bar{e}_k}Z:=\langle DZ,\bar{e}_k\rangle_{\hac_T}$
satisfies the hypotheses of Lemma \ref{l3.1}. Hence, the sequence $([\Lambda_n D^{\bar{e}_k}Z]^{t,x})_{n\in\N}$ converges in $\mathcal{P}_0$ to an element denoted by $\Phi_{t,x}^{(\bar{e}_k)}$. Since the Malliavin derivative is a closed operator, we have $D^{\bar{e}_k}\Phi_{t,x}=\Phi_{t,x}^{(\bar{e}_k)}$, with $\Phi_{t,x}$ given in Lemma \ref{l3.1}, and $\Phi_{t,x}^{(\bar{e}_k)}=\langle \Phi_{t,x}^{(1)},\bar{e}_k\rangle_{\hac_T}$,
with $\Phi_{t,x}^{(1)}$ defined in Lemma \ref{l3.2}.
\smallskip

Let $g$ denote a predictable stochastic process belonging to $L^2(\Omega\times[0,T];\hac)$. Using the stochastic integration theory developed for instance in \cite{dapratozabczyk}, the integral of $g$ with respect to the cylindrical Brownian motion $\{W_t, t\in[0,T]\}$ described in Section \ref{s1} is well-defined, as follows:
\begin{equation}
\label{3.6}
(g\cdot W)_t := \int_0^t g(s) d W_s := \sum_{k\in\N} \int_0^t \langle g(s,\ast),e_k(\ast)\rangle_\hac d W^k_s, \ t\in[0,T],
\end{equation}
where $(e_k)_{k\in\N}$ is a CONS of $\hac$. 

The next proposition provides an extension of \cite[Proposition 2.6]{dalangquer} to the stochastic integral in \cite{conusdalang}. This is only for the class of integrands $\Phi_{t,x}$ defined in Lemma \ref{l3.1} though.

\begin{proposition}
\label{p3.1}
Let $\Lambda$ fulfill {\bf(A3)} and either {\bf(A4)} or {\bf(A5)}. Let $Z$ be a stochastic process satisfying conditions {\bf(A1)} and {\bf(A2)} Fix $(t,x)\in[0,T]\times \Rd$ and consider the stochastic process $\Phi_{t,x}$ defined in Lemma \ref{l3.1}. Then
\begin{equation*}
(\Lambda(t-\cdot,x-\ast)\cdot M^Z)_t = (\Phi_{t,x}\cdot W)_t,\ t\in[0,T].
\end{equation*}
where the expression on the left-hand side refers to the integral of  Conus and Dalang (see \cite[Theorem 3.1]{conusdalang}), while on the right-hand side, it refers to the integral defined in \eqref{3.6}.
\end{proposition}
\noindent{\it Proof.} From Lemma \ref{l3.1}, part 2, we have that $[\Lambda_nZ]^{t,x}\in\mathcal{P}_0\subseteq L^2(\Omega;{\mathcal{H}_T})$. Consequently, the stochastic integral in \eqref{3.6} exists for 
$g:=[\Lambda_nZ]^{t,x}$ and it satisfies the isometry property
\begin{align*}
\E\Big[\big([\Lambda_nZ]^{t,x}\cdot W\big)_t^2\Big] &= \E\Big[\int_0^t \|\Lambda_n(t-s\cdot,x-\ast)Z(s,\ast)\|_\mathcal{H}^2 ds\Big]\\
& = \|[\Lambda_nZ]^{t,x}\|_0^2 = \|\Lambda_n(t-\cdot,x-\ast)\|_{0,Z}^2. 
\end{align*}
Moreover, $(\Lambda_n(t-\cdot,x-\ast)\cdot M^Z)_t$ is well-defined as a Walsh's stochastic integral with respect to the martingale measure $M^Z$. According to \cite[Proposition 2.6(a)]{dalangquer}
\begin{equation*}
(\Lambda_n(t-\cdot,x-\ast)\cdot M^Z)_t = ([\Lambda_nZ]^{t,x}\cdot W)_t.
\end{equation*} We can now  pass to the limit as $n\to\infty$  and notice that $[\Lambda_nZ]^{t,x}$ converges in $L^2(\Omega;\hac_T)$ to $\Phi_{t,x}$. We obtain,
\begin{equation*}
(\Phi_{t,x}\cdot W)_t := L^2(\Omega)-\lim_{n\to\infty}\big([\Lambda_nZ]^{t,x}\cdot W\big)_t. 
\end{equation*}
On the other hand, for the stochastic integral in \cite{conusdalang} we have
\begin{equation*}
(\Lambda(t-\cdot,x-\ast)\cdot M^Z)_t = L^2(\Omega)-\lim_{n\to\infty} (\Lambda_n(t-\cdot,x-\ast)\cdot M^Z)_t .
\end{equation*}
This ends the proof.

\cqd

In the next proposition, we prove the equality between the divergence operator (also called Skorohod integral) applied to the process $\Phi_{t,x}$ and the stochastic integral $(\Phi_{t,x} \cdot W)_t$.

\begin{proposition}
\label{p3.2}
The assumptions are the same as in Proposition \ref{p3.1}. Fix $(t,x)\in[0,T]\times \Rd$. The stochastic process $\Phi_{t,x}$ derived in Lemma \ref{l3.1} part 3 satisfies
\begin{equation}
\label{3.7}
\delta(\Phi_{t,x})= (\Phi_{t,x}\cdot W)_t,
\end{equation}
where $\delta$ denotes the Skorohod integral.
\end{proposition}
\noindent{\it Proof.}  We follow a similar approach as in \cite[Section 1.3.2]{nualart}. Let $g = 1_{(a,b]}1_{A}X$, where $0\leq a<b\leq t$, $A\in\mathcal{B}_b(\Rd)$ and $X$ is a bounded
 and $\mathcal{F}_a$-measurable random variable. Assume first that $X\in\D^{1,2}$. Then \cite[(1.44)]{nualart} yields
 \begin{equation*}
\delta(g) = XF(1_{(a,b]}1_{A}).
\end{equation*}
Since $\D^{1,2}$ is dense in $L^2(\Omega)$ and $\delta$ is closed, this equality extends to $X\in L^2(\Omega)$, which are $\mathcal{F}_a$-measurable. 

On the other hand $(g\cdot W)_t = XF(1_{(a,b]}1_{A})$ as it is shown for instance in \cite[p.\ 11]{dalangquer}. By linearity of the integral operators we see that \eqref{3.7} holds for a suitable class of elementary processes. 

We know that  $[\Lambda_nZ]^{t,x}\in\mathcal{P}_0$. Therefore, there exists a sequence of elementary processes $(g_{n,m}^{t,x})_{m\in\N}$ converging to $[\Lambda_nZ]^{t,x}$ in $L^2(\Omega,{\mathcal{H}_T})$ as $m\to\infty$. Since the operator $\delta$ is closed, we obtain
\begin{equation*}
([\Lambda_nZ]^{t,x}\cdot W)_t = L^2(\Omega)-\lim_{m\to\infty} \left(g_{n,m}^{t,x}\cdot W\right)_t = L^2(\Omega)-\lim_{m\to\infty} \delta(g_{n,m}^{t,x})= \delta([\Lambda_nZ]^{t,x}).
\end{equation*}
Finally, using once again that $\delta$ is closed, we have
\begin{equation*}
(\Phi_{t,x}\cdot W)_t = L^2(\Omega)-\lim_{n\to\infty} \left([\Lambda_nZ]^{t,x}\cdot W\right)_t = L^2(\Omega)-\lim_{n\to\infty} \delta([\Lambda_nZ]^{t,x}) = \delta(\Phi_{t,x}).
\end{equation*}
This finishes the proof of the proposition.

\cqd
%%%%%%%
%%%%%%% Section 4
%%%%%%%
 \section{Malliavin derivatives of stochastic and\break pathwise integrals}
 \label{s4}
In this section, we state conditions for commuting the Malliavin derivative operator with two types of integrals: the class of stochastic integrals studied in Section \ref{s3} and the pathwise integrals of \cite{conusdalang} (see \eqref{1200}). For the former we rely on \cite[Proposition 1.3.2]{nualart} and we check that the assumptions of this proposition are satisfied by the relevant integrands. As for the latter, we give a direct proof.

Throughout the section, we fix $\Lambda$ satisfying the assumption {\bf(A3)} and either {\bf(A4)} or {\bf(A5)} and a stochastic process $Z$ satisfying {\bf(A1)} and {\bf(A2)}. For any $(t,x)\in[0,T]\times \Rd$, we shall consider the stochastic process $\Phi_{t,x}$ given in Lemma \ref{l3.1}. We will use the notation 
\begin{equation*}
\int_0^t\int_\Rd \Lambda(t-s,x-z)Z(s,z)M(ds,dz)
\end{equation*}
to refer to each of the stochastic integrals $\delta(\Phi_{t,x})$, $(\Phi_{t,x}\cdot W)_t$, $(\Lambda(t-\cdot,x-\ast)\cdot M^Z)_t$ considered in Section \ref{s3}.
In fact, owing to Propositions \ref{p3.1}, \ref{p3.2} they coincide.

 \begin{proposition}
\label{p4.1}
We assume that $\Lambda$ satisfies the assumption {\bf(A3)} and either {\bf(A4)} or {\bf(A5)}. Consider a stochastic process $Z$ satisfying {\bf(A1)} and {\bf(A2)} and such that for any $(t,x)\in[0,T]\times \Rd$, $Z(t,x)\in \D^{1,2}$. Suppose also that $DZ$ fulfills the assumptions  {\bf(A8)} and {\bf(A9)} with $\caA={\mathcal{H}_T}$. Then, for every $(t,x)\in[0,T]\times \Rd$,
\begin{equation*}
\int_0^t\int_{\Rd} \Lambda(t-s,x-z)Z(s,z)M(ds,dz) \in \D^{1,2}
\end{equation*}
  and
  \begin{align}
  \label{4.1}
  &D\left(\int_0^t\int_\Rd \Lambda(t-s,x-z)Z(s,z)M(ds,dz)\right)\nonumber\\
    & \quad \quad \quad = \Lambda(t-\cdot,x-\ast)Z(\cdot,\ast) + \int_0^t\int_\Rd \Lambda(t-s,x-z)DZ(s,z)M(ds,dz),
  \end{align}  
  where the integral in the right-hand side of \eqref{4.1} is the Hilbert space valued stochastic integral $\Lambda(t-\cdot,x-\ast)\cdot M^{DZ}$ given in Theorem \ref{t2.1}.
  \end{proposition}
  
\noindent{\it Proof.} From Lemma \ref{l3.2} it follows that $\Phi_{t,x}\in \D^{1,2}({\mathcal{H}_T})$. Fix a CONS of $\hac_T$ that we denote by $(\bar{e}_k)_{k\in\N}$. As has been pointed out in Section \ref{s3}, the real-valued process $D^{\bar{e}_k}\Phi_{t,x}=\langle D\Phi_{t,x}, \bar{e}_k\rangle_{\mathcal{H}_T}$ belongs to $\mathcal{P}_0$. Moreover, by the results of that section, it also belongs to the domain of the divergence operator. Thus, the assumptions of \cite[Proposition 1.3.2]{nualart} are fulfilled and hence we have
\begin{align*}
D^{\bar{e}_k}\bigg(\int_0^t\int_\Rd \Phi_{t,x}(s,z)M(ds,dz)\bigg)&= 
D^{\bar{e}_k}(\delta(\Phi_{t,x}))\\
&= \langle\Phi_{t,x}, \bar{e}_k\rangle_{\mathcal{H}_T} + \delta(D^{\bar{e}_k}\Phi_{t,x}),
\end{align*}
for any $k\in\N$. This proves \eqref{4.1}.

\cqd

Our next aim is to prove a result on commutation of the Malliavin derivative operator with the pathwise integral \eqref{1200}. In \cite[Lemma 2.2]{stroockkusuoka} a similar question is analyzed. 
However, that version seems not to be directly applicable to our context.

\begin{proposition}
\label{p4.2}
  Let $\Lambda$ fullfil {\bf(A6)} and either {\bf(A4)} or {\bf(A7)}. Let $Z$ be a stochastic process satisfying the same assumptions as in Proposition \ref{p4.1}. Then for all
  $(t,x)\in[0,T]\times \Rd$, 
  \begin{equation*}
   \int_0^t\int_\Rd \Lambda(t-s,x-z) Z(s,z) dzds \in \D^{1,2} 
   \end{equation*}
  and 
  \begin{equation}
  \label{4.2}
  D\bigg(\int_0^t\int_\Rd \Lambda(t-s,x-z) Z(s,z) dzds\bigg) = \int_0^t\int_\Rd \Lambda(t-s,x-z)DZ(s,z) dzds.
  \end{equation}
  \end{proposition}
\noindent{\it Proof.}
Let $\Lambda=1_{t-(a,b]}1_{\{x\}-A}$ for some $0\leq a<b\leq t$ and $A\in\mathcal{B}_b(\Rd)$ where $\{x\}-A=\{x-z, z\in A\}$. 
In this case, formula \eqref{4.2} reads
\begin{equation}
\label{4.3}
	D\bigg(\int_a^b\int_A Z(s,z) dz ds\bigg) = \int_a^b\int_A DZ(s,z) dzds
	\end{equation}  
almost surely. This follows from the arguments in the above-mentioned reference \cite{stroockkusuoka}. We notice that a direct proof of \eqref{4.3} can also be done using the definition of the Malliavin operator as a directional derivative.

In the next step, we consider  $\Lambda\in L^2([0,T];L^1(\Rd))$, i.e.\ $\int_0^T(\int_\Rd |\Lambda(t,x)|d x)^2d t < \infty$. We recall that according to \cite[(3.13)]{conusdalang}, the pathwise integral is almost surely well defined as a pathwise Lebesgue integral. Linear combinations of products of indicator functions as those considered in the previous step, are dense in 
$L^2([0,T];L^1(\Rd))$. Let $(\Lambda_n)_{n\in\N}$ be a sequence of such simple functions converging to $\Lambda$ in $L^2([0,T];L^1(\Rd))$. Then,
\begin{align*}
  & \E\left[\left(\int_0^t\int_\Rd \left(\Lambda(t-s,x-z)-\Lambda_n(t-s,x-z)\right) Z(s,z) dzds\right)^2\right] \\
  & \leq T\E\Big[\int_0^t\int_\Rd\int_\Rd \left(\Lambda(t-s,x-z)-\Lambda_n(t-s,x-z)\right)\\
  &\qquad \times \left(\Lambda(t-s,x-y)-\Lambda_n(t-s,x-y)\right) Z(s,z)  Z(s,y) dydzds\Big] \\
  & \leq C\sup_{(t,x)\in[0,T]\times \Rd} \E(|Z(t,x)|^2)\\
  &\qquad\times \int_0^t\left(\int_\Rd |\Lambda(t-s,x-z)-\Lambda_n(t-s,x-z)|\ dz\right)^2ds,
\end{align*}
which goes to zero as $n\to\infty$. 

Using similar arguments,
\begin{align*}
  & \E\left[\left(\int_0^t\int_\Rd \Lambda(t-s,x-z)D Z(s,z) dzds\right.\right.\\
  &\left.\left.\qquad - D\int_0^t\int_\Rd \Lambda_n(t-s,x-z) Z(s,z) dzds\right)^2\right] \\
  & = \E\left[\left(\int_0^t\int_\Rd \big(\Lambda(t-s,x-z)-\Lambda_n(t-s,x-z)\big)D Z(s,z) dzds\right)^2\right] \\
  & \leq C \sup_{(t,x)\in[0,T]\times \Rd} \E\left[\|DZ(t,x)\|^2_{\hac_T}\right] \\
  &\qquad\times \int_0^t\left(\int_\Rd |\Lambda(t-s,x-z)-\Lambda_n(t-s,x-z)| dz\right)^2ds,
\end{align*}
where in the first equality we have used the first step of this proof. The last term goes to zero as $n\to\infty$. Since $D$ is a closed operator, the Proposition holds
for $\Lambda\in L^2([0,T];L^1(\Rd))$.

Finally, assume that $\Lambda$ satisfies the assumptions of the Proposition. Let $\Lambda_n\in L^2([0,T],L^1(\Rd))$, $n\in\N$ be as in \eqref{eq:Lambda_n}. 
Then, according to \cite{conusdalang} (see also Section \ref{s1}) we have,
\begin{align*}
  & \E\left[\left(\int_0^t\int_\Rd \big(\Lambda(t-s,x-z)-\Lambda_n(t-s,x-z)\big) Z(s,z) dzds\right)^2\right]\\
  & = \int_0^t\int_\Rd |\tf\Lambda(t-s)(\eta)|^2|\tf\zeta_n(\eta)-1|^2 \nu^{Z}_s(d\eta)ds.
\end{align*}
This goes to zero as $n\to\infty$, by dominated convergence. 

Similarly,
\begin{align*}
  & \E\left[\left(\int_0^t\int_\Rd \Lambda(t-s,x-z)DZ(s,z) dzds\right.\right.\\
  &\left.\left.\qquad - D\int_0^t\int_\Rd \Lambda_n(t-s,x-z) Z(s,z) dzds\right)^2\right] \\
  & = \E\left[\left(\int_0^t\int_\Rd \big(\Lambda(t-s,x-z)-\Lambda_n(t-s,x-z)\big)DZ(s,z) dzds\right)^2\right] \\
  & = \int_0^t\int_\Rd |\tf\Lambda(t-s)(\eta)|^2|\tf\zeta_n(\eta)-1|^2 \nu^{DZ}_s(d\eta)ds,
\end{align*}
converges to zero as $n\to\infty$. Notice that by Theorem \ref{t2.2}, the integrals involved in these computations exist. By the closedness of the Malliavin derivative operator, 
we conclude the proof.

\cqd
%%%%%%%%
%%%%%%%%%MALLIAVIN DERIVATIVE OF THE SOLUTION
%%%%%%%%
%%%%%%%%

\section{Malliavin differentiability of the solution of the SPDE}
\label{s5}

This section is devoted to prove that the solution to the stochastic partial differential equation  \eqref{i.1} at a given point $(t,x)\in[0,T]\times \Rd$ is differentiable in Malliavin's sense. We also derive an SPDE satisfied by the $\hac_T$--valued stochastic process $\{Du(t,x), (t,x)\in[0,T]\times \Rd\}$. This general result applies in particular to the solution of the stochastic wave equation in any spatial dimension.

%%%%%%%Remark on the equations to be considered
%%%%%%%%
%%%%%%%%
It is assumed that $G$ satisfies {\bf (A3)}, {\bf (A6)} and either {\bf (A4)} or  {\bf (A5)} and {\bf (A7)}.
For its further use, we introduce an SPDE more general than \eqref{i.1}, as follows. Let $h\in\hac_T$ and consider
\begin{align}
\label{5.a}
u^h(t,x)&=\int_0^t\int_{\Rd} G(t-s,x-z)\sigma(u(s,z)) M(ds,dz)\nonumber\\
&+\int_0^t\langle G(t-s,x-\ast)\sigma(u^h(s,\ast)),h(s)\rangle_{\hac} ds\nonumber\\
&+\int_0^t\int_{\Rd} G(t-s,x-z) b(u(s,z)) dz ds.
\end{align}
It is easy to check that the Picard iterations $\{u^{m,h}(t,x), (t,x)\in[0,T]\times \Rd\}$, $m\in\N$, satisfy the $S$ property of \cite[Lemma 4.5]{conusdalang}. With this,
an easy extension of \cite[Theorem 4.2, Theorem 4.8]{conusdalang} provides existence (and uniqueness) of a random field solution $\{u^h(t,x), (t,x)\in[0,T]\times\Rd\}$
to \eqref{5.a}. Moreover,
\begin{equation}
\label{5.b}
\sup_{(t,x)\in[0,T]\times\Rd}\sup_{\Vert h\Vert_{\hac_T}\le c} E\left[\left\vert u^h(t,x)\right\vert^2\right] < \infty.
\end{equation}
The details of the proof are left to the reader.

Owing to the results proved in Section \ref{s3}, the stochastic integral in \eqref{5.a} can be interpreted either as a Conus-Dalang's integral, a Skorohod integral, or as
\begin{equation*}
\int_0^t\int_{\Rd} G(t-s,x-z)\sigma(u(s,z)) M(ds,dz):=\sum_{k\in\N}\int_0^t \langle G(t-s,x-\ast)\sigma(u^h(s,\ast)),e_k\rangle_{\hac} dW_s^k,
\end{equation*}
with $(e_k)_{k\in\N}$ a CONS of $\hac$ and $(W^k_t, t\in[0,T])_{k\in\N}$ a sequence of independent standard Brownian motions.

Similarly
\begin{equation*}
\int_0^t\langle G(t-s,x-\ast)\sigma(u^h(s,\ast)),h(s)\rangle_{\hac} ds
=\sum_{k\in\N}\int_0^t \langle G(t-s,x-\ast)\sigma(u^h(s,\ast)),e_k\rangle_{\hac} h^k(s) ds,
\end{equation*}
where $h^k(s)=\langle h(s),e_k\rangle_{\hac}$, $k\in\N$.

Throughout the section, we shall use the abstract Wiener space $(\Omega,\mathbb{H}, \mathbb{P})$ and the isometry between the spaces $\hac_T$ and $\mathbb{H}$ defined in Section \ref{s1}. 
%%%%%%% Back to the proof of Theorem \ref{t5.1}

The objective is to prove the following.

\begin{theorem}
\label{t5.1} We assume that $G$ satisfies the assumptions {\bf(A3)}, {\bf(A6)} and either {\bf(A4)} or {\bf(A5)} and {\bf(A7)}. We also suppose that the coefficients $\sigma$ and $b$ are continuously differentiable real-valued functions with bounded derivatives. Then for any $(t,x)\in[0,T]\times \Rd$, $u(t,x)\in\mathbb{D}^{1,2}$. Moreover, the stochastic process 
$\{Du(t,x), (t,x)\in[0,T]\times \Rd\}$ satisfies the SPDE
\begin{align}
\label{5.1}
Du(t,x)&= G(t-\cdot,x-\ast)\sigma(u(\cdot,\ast))\nonumber\\
&+\int_0^t\int_{\Rd} G(t-s,x-z)\sigma^\prime(u(s,z)) Du(s,z) M(ds,dz)\nonumber\\
&+\int_0^t\int_{\Rd} G(t-s,x-z)b^\prime(u(s,z)) Du(s,z)\ dz ds,
\end{align}
where $G(t-\cdot,x-\ast)\sigma(u(\cdot,\ast))$ is the stochastic process derived in Lemma \ref{l3.1} Section \ref{s3}, for $\Lambda:=G$ and $Z:=\sigma(u)$.
\end{theorem}

The proof of Theorem \ref{t5.1} will be carried out in two steps. Firstly, we will show that $u(t,x)\in\mathbb{D}^{1,2}$ and in a second step, we shall establish \eqref{5.1}. The proof of the former statement relies on \cite[Lemma 1.2.3]{nualart}. For the sake of completeness, we quote this result.
\begin{lemma}
\label{l5.1}
  Let $(F_n)_{n\in\N}$ be a sequence in $\D^{1,2}$ such that $\lim_{n\to\infty}F_n= F$ in $L^2(\Omega)$  and $\sup_{n\in\N} \E[\|DF_n\|^2_{\mathcal{H}_T}]<\infty$. Then $F\in\D^{1,2}$ and the sequence $(DF_n)_{n\in\N}$ converges to $DF$ in the weak topology of $L^2(\Omega;{\mathcal{H}_T})$.
\end{lemma}

This Lemma will be applied to the sequence $F_n:=u_n(t,x)$, $n\in\N$, where $(t,x)\in[0,T]\times \Rd$ is fixed, and $u_n(t,x)$ is given by the solution to the evolution equation
\begin{align}
\label{5.2}
  u_n(t,x) &= \int_0^t\int_{\Rd} G_n(t-s,x-z)\sigma(u_n(s,z))M(ds,dz) \notag\\
  					& + \int_0^t\int_{\Rd} G_n(t-s,x-z)b(u_n(s,z))dz ds, 
\end{align}
with $G_n$ defined as in \eqref{eq:Lambda_n}.

Assume that the functions $\sigma$, $b$ are Lipschitz continuous. Since $G_n(t)\in\mathcal{S}(\Rd)$, the stochastic integral in \eqref{5.2} is a Walsh's integral (see \cite{walsh}). It is well-known that \eqref{5.2} has a unique random field solution, and that it satisfies the $S$--property. In particular for each $n\in\N$, the process $\{Z(t,x):=u_n(t,x), (t,x)\in[0,T]\times\Rd\}$ satisfies the assumptions {\bf(A1)}, {\bf(A2)}. For a proof of these results, we can proceed as in \cite[Theorem 13]{dalang}.
 \smallskip

From the proof of Proposition 7.1 in \cite{sanzbook} we obtain the following.
\begin{proposition}
\label{p5.1}
Let $G_n:=G\ast\zeta_n$ be as in \eqref{eq:Lambda_n}. Assume that the coefficients $\sigma$, $b$ in \eqref{5.2} are continuously differentiable with bounded derivatives. Then for each $n\in\N$ and every $(t,x)\in[0,T]\times\Rd$, the 
random variable $u_n(t,x)$ belongs to $\mathbb{D}^{1,2}$. Moreover, the $\hac_T$-valued stochastic process $\{Du_n(t,x), (t,x)\in[0,T]\times\Rd\}$ is the solution to the SPDE
\begin{align}
\label{5.3}
	Du_n(t,x) = & G_n(t-\cdot,x-\ast)\sigma(u_{n}(\cdot,\ast)) \notag\\
              & + \int_0^t\int_{\Rd} G_n(t-s,x-z)\sigma'(u_n(s,z))Du_n(s,z)M(ds,dz) \notag\\
              & + \int_0^t\int_{\Rd} G_n(t-s,x-z)b'(u_n(s,z))Du_n(s,z) dz ds.
\end{align}
\end{proposition}
Next, we study the convergence of the sequence of processes $(u_n)_{n\in\N}$ to $u$.

\begin{proposition}
\label{p5.2}
We assume that $G$ satisfies the  hypotheses of Theorem \ref{t5.1}. Moreover, we suppose that the functions $\sigma$ and $b$ are Lipschitz continuous. Then we have
\begin{equation*}
 \lim_{n\to\infty}\sup_{(t,x)\in[0,T]\times\Rd}\E\big[|u_n(t,x)-u(t,x)|^2\big] = 0. 
 \end{equation*}
 \end{proposition}
\noindent{\it Proof.} We start by proving that
\begin{equation}
\label{5.4}
  \sup_{n\in\N} \sup_{(t,x)\in[0,T]\times\Rd}\ \E\big[|u_n(t,x)|^2\big] < \infty.
\end{equation}
Indeed, from \eqref{5.2} it follows that $\E\big[|u_n(t,x)|^2\big]\leq 2(I_{1,n}(t,x) + I_{2,n}(t,x))$, for every $(t,x)\in[0,T]\times \Rd$, where 
\begin{align*}
	I_{1,n}(t,x)  & = \E\Bigg[\bigg(\int_0^t\int_{\Rd} G_n(t-s,x-z)\sigma(u_n(s,z)) M(ds,dz)\bigg)^2\Bigg]
\intertext{and}
	I_{2,n}(t,x)  & = \E\Bigg[\bigg(\int_0^t\int_\Rd G_n(t-s,x-z)b(u_n(s,z)) dz ds\bigg)^2\Bigg].
\end{align*}
Notice that the inequalities \eqref{eq:boundCDI}, \eqref{1.3} also hold with $\Lambda$ replaced by $G_n(t-\cdot,x-\ast)$. Then, by taking $Z(t,x):=\sigma(u_n(t,x))$ and $Z(t,x):=b(u_n(t,x))$, respectively, we obtain
\begin{align*}
  I_{1,n}(t,x)  & \leq \int_0^t\sup_{(r,y)\in[0,s]\times\Rd} \E[\sigma(u_n(r,y))^2]
  \sup_{\eta\in\Rd}\int_{\Rd} |\tf G_n(t-s)(\xi+\eta)|^2\mu(d\xi) ds\\
                & \leq C\int_0^t\sup_{(r,y)\in[0,s]\times\Rd} \E[(1+u_n(r,y))^2]
                \sup_{\eta\in\Rd}\int_{\Rd} |\tf G(t-s)(\xi+\eta)|^2\mu(d\xi) ds\\
                & \leq C\int_0^t \bigg(1+ \sup_{(r,y)\in[0,s]\times\Rd} \E[u_n(r,y)^2]\bigg)J_1(t-s) ds,
\end{align*}
and 
\begin{align*}
  I_{2,n}(t,x)  & \leq \int_0^t \sup_{(r,y)\in[0,s]\times\Rd} \E[b(u_n(r,y))^2] \sup_{\eta\in\R} |\tf G_n(t-s)(\eta)|^2 ds\\
                & \leq C \int_0^t  \bigg(1+ \sup_{(r,y)\in[0,s]\times\Rd} \E[u_n(r,y)^2]\bigg) J_2(t-s) ds,
\end{align*}
where the functions $J_1$ and $J_2$ are defined in \eqref{eq:J1} and \eqref{eq:J2}, respectively with $\Lambda$ replaced by $G$. 
 This yields
 \begin{align*}
 \sup_{(r,y)\in[0,t]\times\Rd} \E[|u_n(r,y)|^2]& \leq C\int_0^t \bigg(1+ \sup_{(r,y)\in[0,s]\times\Rd} \E[u_n(r,y)^2]\bigg)\\
 &\quad \times\big(J_1(t-s)+J_2(t-s)\big) ds. 
 \end{align*}
Using the version of Gronwall's Lemma in \cite[Lemma 15]{dalang} along with \eqref{eq:2.17'''} yields \eqref{5.4}.

Next, we show the assertion of the proposition. Using equations \eqref{i.1} and \eqref{5.2}, we have
\[ \E\big[|u(t,x)-u_n(t,x)|^2\big] \leq C(T_{1,n}(t,x)+T_{2,n}(t,x)+T_{3,n}(t,x)+T_{4,n}(t,x)), \]
where
\begin{align*}
	T_{1,n}(t,x) & = \E\Bigg[\bigg(\int_0^t\int_\Rd G_n(t-s,x-z)\\
	&\qquad\times\big(\sigma(u_n(s,z))-\sigma(u(s,z))\big)M(ds,dz)\bigg)^2\Bigg], \\
	T_{2,n}(t,x) & = \E\Bigg[\bigg(\int_0^t\int_\Rd \big(G_n(t-s,x-z)-G(t-s,x-z)\big)\\
	&\qquad \times\sigma(u(s,z))M(ds,dz)\bigg)^2\Bigg], \\
	T_{3,n}(t,x) & = \E\Bigg[\bigg(\int_0^t\int_\Rd G_n(t-s,x-z)\big(b(u_n(s,z))-b(u(s,z))\big) dz ds\bigg)^2\Bigg],\\
	T_{4,n}(t,x) & = \E\Bigg[\bigg(\int_0^t\int_\Rd \big(G_n(t-s,x-z)-G(t-s,x-z)\big)\\
	&\qquad\times b(u(s,z)) dz ds\bigg)^2\Bigg].	
\end{align*}
For the terms $T_{1,n}(t,x)$, $T_{2,n}(t,x)$, we apply the inequality \eqref{eq:boundCDI} in the following situations. For the former term, we replace $\Lambda$ by
 $G_n(t-\cdot,x-\ast)$ and take $Z:=\sigma(u_n)-\sigma(u)$; for the latter, we replace $\Lambda$ by $[G_n-G](t-\cdot,x-\ast)$ and take $Z:=\sigma(u)$. This yields
\begin{align*}
  & T_{1,n}(t,x)   
   \leq \int_0^t\sup_{(r,y)\in[0,s]\times\Rd} \E\Big[\big|\sigma(u_n(r,y))-\sigma(u(r,y))\big|^2\Big]\\
  &\qquad\times \sup_{\eta\in\Rd}\int_{\Rd} |\tf G_n(t-s)(\xi+\eta)|^2\mu(d\xi) ds \\
  & \leq C\int_0^t\sup_{(r,y)\in[0,s]\times\Rd} \E\Big[\big|u_n(r,y)-u(r,y)\big|^2\Big]J_1(t-s) ds,
\end{align*}
and
\begin{align*}
  T_{2,n}(t,x)	& = \int_0^t\int_\Rd \big|\tf G_n(t-s)(\xi)-\tf G(t-s)(\xi)\big|^2\mu^{\sigma(u)}_s(d\xi) ds \\
                & = \int_0^t\int_\Rd \big|\tf G(t-s)(\xi)\big|^2\big|\tf\zeta_n(\xi)-1\big|^2\mu^{\sigma(u)}_s(d\xi) ds.
\end{align*}                
For the term $T_{3,n}(t,x)$, we apply \eqref{1.3} with $\Lambda$ replaced by $G_n(t-\cdot,x-\ast)$ and $Z:=b(u_n)-b(u)$. For $T_{4,n}(t,x)$, we proceed similarly with $\Lambda$ replaced by $[G_n-G](t-\cdot,x-\ast)$ and $Z:=b(u)$,
respectively. We obtain
\begin{align*}
  T_{3,n}(t,x)& \leq \int_0^t \sup_{(r,y)\in[0,s]\times\Rd} \E\Big[\big|b(u_n(r,y))-b(u(r,y))\big|^2\Big]
  \sup_{\eta\in\Rd} |\tf G_n(t-s)(\eta)|^2 ds\\
                & \leq C\int_0^t \sup_{(r,y)\in[0,s]\times\Rd} \E\Big[\big|u_n(r,y)-u(r,y)\big|^2\Big] J_2(t-s) ds,\\
  T_{4,n}(t,x) 	& = \int_0^t\int_\Rd \big|\tf G_n(t-s)(\eta)-\tf G(t-s)(\eta)\big|^2\nu^{\sigma(u)}_s(d\eta) ds\\
                & = \int_0^t\int_\Rd |\tf G(t-s)(\eta)|^2|\tf\zeta_n(\eta)-1|^2\nu^{\sigma(u)}_s(d\eta) ds.
\end{align*}
The terms $T_{2,n}(t,x)$, $T_{4,n}(t,x)$ converge to zero as $n\to\infty$ uniformly in $(t,x)\in[0,T]\times\Rd$, by dominated convergence. Hence, altogether we have
\begin{align*}
  &\sup_{(r,y)\in[0,t]\times\Rd} \E\big[|u(r,y)-u_n(r,y)|^2\big] \\
            & \quad\leq C_n + C\int_0^t \sup_{(r,y)\in[0,s]\times\Rd} \E\big[|u(r,y)-u_n(r,y)|^2\big](J_1(t-s)+J_2(t-s)) ds,
\end{align*}
where $C_n$ tends to $0$ as $n\to\infty$ uniformly in $(t,x)\in[0,T]\times\Rd$. An application of Gronwall's Lemma yields the assertion.

\cqd
%%%%%%%%

The next proposition provides the last ingredient for the application of Lemma \ref{l5.1}.

\begin{proposition}
\label{p5.3}
  With the same assumptions as in Theorem \ref{t5.1}, we have
  \begin{equation*}
  \sup_{n\in\N} \sup_{(t,x)\in[0,T]\times\Rd} \E\big[\|Du_n(t,x)\|^2_{\mathcal{H}_T}\big]<\infty. 
  \end{equation*}
  \end{proposition}
\noindent{\it Proof.}
Fix $(t,x)\in[0,T]\times\Rd$. We bound the $L^2(\Omega;{\mathcal{H}_T})$-norm of each term on the right-hand side of \eqref{5.3}. For the first term, we apply \eqref{3.1} with $\phi:=G_n(t-\cdot,x-\ast)$ and $Z:=\sigma(u_n)$ and then, \eqref{eq:normP0Z}, \eqref{eq:boundCDI} with $g=\Lambda:=G_n(t-\cdot,x-\ast)$.  By the properties of $\sigma$, we obtain
\begin{align*}
  & \E\big[\|G_n(t-\cdot,x-\ast)\sigma(u_n(\cdot,\ast))\|^2_{\hac_T}\big]\\
  & = \int_0^t\int_\Rd |\tf G_n(t-s)(\xi)|^2\mu_s^{\sigma(u_n)}(d\xi)ds\\
  & \leq \int_0^t\E[\sigma(u_n(s,0))^2]\sup_{\eta\in\Rd}\int_\Rd |\tf G_n(t-s)(\xi+\eta)|^2\mu(d\xi)ds\\
  & \leq C\bigg(1+\sup_{(r,y)\in[0,T]\times\Rd}\E\big[u_n(r,y)^2\big]\bigg)\int_0^t\sup_{\eta\in\Rd}\int_\Rd|\tf G(t-s)(\xi+\eta)|^2\mu(d\xi)ds.
\end{align*}
The last term is uniformly bounded in $n\in\N$ and $(t,x)\in[0,T]\times\Rd$, due to \eqref{5.4} and assumption {\bf(A3)}.

For the second term on the right-hand side of \eqref{5.3}, we apply \eqref{2.1} with $\Lambda$ replaced by $G_n(t-\cdot,x-\ast)$, $Z:=\sigma^\prime(u_n) Du_n$ and $\mathcal{A}=\hac_T$.
Since $\sigma^\prime$ is bounded, we obtain
\begin{align*}
  \E\Bigg[\bigg\| & \int_0^t\int_\Rd G_n(t-s,x-z)\sigma'(u_n(s,z))Du_n(s,z)M(ds,dz)\bigg\|^2_{\mathcal{H}_T}\Bigg] \\
  & \leq \int_0^t \E\big[\|\sigma'(u_n(s,0))Du_n(s,0)\|^2_{\mathcal{H}_T}\big]\sup_{\eta\in\Rd}\int_\Rd|\tf G_n(t-s)(\xi+\eta)|^2\mu(d\xi)ds\\
  & \leq C\int_0^t \sup_{(r,y)\in[0,s]\times\Rd}\E\big[\|Du_n(r,y)\|^2_{\mathcal{H}_T}\big]J_1(t-s) ds.
\end{align*}
Finally, applying \eqref{2.3} with $\Lambda$ replaced by $G_n(t-\cdot,x-\ast)$, $Z:=b^\prime(u_n) Du_n$ and $\mathcal{A}=\hac_T$ yields
\begin{align*}
  \E\Bigg[\bigg\| & \int_0^t\int_\Rd G_n(t-s,x-z)b'(u_n(s,z))Du_n(s,z) dz ds\bigg\|^2_{\mathcal{H}_T}\Bigg] \\
  & \leq \int_0^t \E\big[\|b'(u_n(s,0))Du_n(s,0)\|^2_{\mathcal{H}_T}\big]\sup_{\eta\in\Rd}|\tf G_n(t-s)(\eta)|^2ds\\
  & \leq C\int_0^t \sup_{(r,y)\in[0,s]\times\Rd}\E\big[\|Du_n(r,y)\|^2_{\mathcal{H}_T}\big]J_2(t-s)ds.
\end{align*}
Thus, 
\begin{align*}
&\E\big[\|Du_n(t,x)\|^2_{\mathcal{H}_T}\big]\\
&\qquad \leq C\left[1+\int_0^t\sup_{(r,y)\in[0,s]\times\Rd}\E\big[\|Du_n(r,y)\|^2_{\mathcal{H}_T}\big](J_1(t-s)+J_2(t-s)) ds\right].
\end{align*}
An application of Gronwall's Lemma finishes the proof.

\cqd

Propositions \ref{p5.2}, \ref{p5.3}, along with Lemma \ref{l5.1} yields that $u(t,x)\in\mathbb{D}^{1,2}$ for any $(t,x)\in[0,T]\times \R^d$. This is the first assertion of
Theorem \ref{t5.1}.
\medskip

The rest of this section is devoted to prove that the Malliavin derivative of the process $\{u(t,x), (t,x)\in[0,T]\times \Rd\}$ satisfies \eqref{5.1}. For this, we consider the equation \eqref{i.1} satisfied by this 
process and apply the Malliavin derivative operator to each term. We obtain
\begin{align}
\label{5.5}
Du(t,x)&= D\left(\int_0^t\int_{\Rd} G(t-s,x-z) \sigma(u(s,z)) M(ds,dz)\right)\nonumber\\
&+D\left(\int_0^t\int_{\Rd} G(t-s,x-z) b(u(s,z)) dz ds\right).
\end{align}
Then \eqref{5.1} will follow by applying Propositions \ref{p4.1}, \ref{p4.2}. The rest of this section is devoted to check that the stochastic processes $Z(t,x):=\sigma(u(t,x))$ and
 $Z(t,x):=b(u(t,x))$, $(t,x)\in[0,T]\times\Rd$, satisfy the assumptions of these propositions, respectively.
 \medskip
 
 \begin{lemma}
 \label{l5.2}
 Let $B:\R\rightarrow \R$ be a Lipschitz continuous function. Then the stochastic process $B(u)=\{B(u(t,x)), (t,x)\in[0,T]\times \Rd\}$, where $u=\{u(t,x), (t,x)\in[0,T]\times \Rd\}$ is the solution of \eqref{i.1}, satisfies the assumptions {\bf(A1)}, {\bf(A2)}.
 \end{lemma}
 \noindent{\it Proof.} Since the process $u$ is predictable and $B$ is continuous, $B(u)$ is clearly predictable. The function $B$ has linear growth; along with \eqref{5.b}, this  yields
 \begin{equation*}
 \sup_{(t,x)\in[0,T]\times\Rd}\E\left[B(u(t,x))^2\right]\le C\left[1+\sup_{(t,x)\in[0,T]\times\Rd}\E\left(|u(t,x)|^2\right)\right]<\infty.
 \end{equation*}

The proof of {\bf(A2)} follows from the $S$-property of the process $u$ (see \cite[Definition 4.4, Lemma 4.5 and Theorem 4.2]{conusdalang}).

\cqd
\medskip

\begin{lemma}
\label{l5.3}
Let $B(u)=\{B(u(t,x)), (t,x)\in[0,T]\times \Rd\}$ be as in Lemma \ref{l5.2}. Assume in addition that $B$ is continuous differentiable with bounded derivative. Then the $\hac_T$-valued stochastic process $D(B(u)):\{D(B(u(t,x))), (t,x)\in[0,T]\times\Rd\}$ satisfies {\bf (A8)}, {\bf (A9)}.
\end{lemma}

\noindent{\it Proof.} First, we note that by the construction of the Malliavin derivative based on smooth functionals (see for instance \cite[(1.29)]{nualart}), the stochastic process $D(B(u))$ inherits the predictability property of the process $u$. 
We also notice that by the chain rule of Malliavin calculus, $B(u(t,x))\in \D^{1,2}$ and $D\left(B(u(t,x))\right)=B^\prime(u((t,x)) Du(t,x)$, for any $(t,x)\in[0,T]\times \Rd$. 

We are assuming that $B^\prime$ is bounded. Thus,
\begin{align*}
\E\left[\Vert D\left(B(u(t,x))\right)\Vert_{\hac_T}^2\right]&\le C\E\left[\Vert Du(t,x)\Vert_{\hac_T}^2\right]
\le C \liminf_{n\to\infty}\E\left[\Vert Du_n(t,x)\Vert_{\hac_T}^2\right]\\
&\le C\sup_{n\in\N}\E\left[\Vert Du_n(t,x)\Vert_{\hac_T}^2\right],
\end{align*}
where $u_n(t,x)$ is defined by \eqref{5.2}.  In the second inequality above, we have used that the sequence $(Du_n(t,x))_{n\in\N}$ converges weakly in $\hac_T$ to $Du(t,x)$ along with \cite[Theorem 5, Chapter 10]{lax}.
From Proposition \ref{p5.3}, we conclude
\begin{equation*}
\sup_{(t,x)\in[0,T]\times\Rd}\E\left[\Vert D\left(B(u(t,x))\right)\Vert_{\hac_T}^2\right] <\infty.
\end{equation*}
Hence the stochastic process $D(B(u))$ satisfies {\bf(A8)}.
\smallskip 

Consider the Picard iterations of the processes $\{u(t,x), (t,x)\in[0,T]\times \Rd\}$,  $\{u^h(t,x), (t,x)\in[0,T]\times \Rd\}$, $h\in \hac_T$, that we denote by 
$\{u^m(t,x), (t,x)\in[0,T]\times \Rd\}$,  $\{u^{m,h}(t,x), (t,x)\in[0,T]\times \Rd\}$, $m\ge 1$, respectively. We have the following:
\smallskip

\noindent  {\bf(SP)} for any $m\ge 1$,
the process
\begin{equation*}
\left(u^m(t,x), u^{m,h}(t,x), u^{m-1}(t,x), u^{m-1,h}(t,x), (t,x)\in[0,T]\times\Rd\right),
\end{equation*}
satisfies the $S$-property defined in \cite[Definition 4.4]{conusdalang}. 

Indeed, this can be proved by a recursive argument on $m\ge 1$, following similar arguments as in \cite[Lemma 4.5]{conusdalang}.

Property {\bf(SP)} implies that the process $D(B(u))$ satisfies {\bf (A9)}. 
Indeed, let $(\bar{e}_k)_{k\in\N}$ be a CONS of $\hac_T$. The Malliavin derivative $D^{\bar{e}_k}u(t,x)$ can be obtained as
\begin{equation*}
L^2(\Omega)-\lim_{\ep\to 0}\frac{u^{\ep \bar{e}_k}(t,x)-u(t,x)}{\ep}.
\end{equation*}
Then, using the chain rule of Malliavin calculus and dominated convergence twice, we conclude
\begin{align*}
\label{5.6}
&\E\left[D^{\bar{e}_k}(B(u(t,x))) D^{\bar{e}_k}(B(u(t,x+y)))\right]\nonumber\\
& \quad =\lim_{\ep\to 0}\E\left[B^\prime(u(t,x))\frac{u^{\ep \bar{e}_k}(t,x)-u(t,x)}{\ep}\right.\\
&\left.\quad\quad \quad \times B^\prime(u(t,x+y))\frac{u^{\ep \bar{e}_k}(t,x+y)-u(t,x+y)}{\ep}\right]\nonumber\\
&\quad=\E\left[D^{\bar{e}_k}(B(u(t,0))) D^{\bar{e}_k}(B(u(t,y)))\right],
\end{align*}
where the last equality is a consequence of {\bf(SP)}.
\cqd

%%%%%Section 6 Existence of Density

\section{Existence of density}
\label{s6} 

In this section we consider the solution to the SPDE \eqref{i.1} at a fixed point $(t,x)\in]0,T]\times \Rd$ in the particular
case where $\sigma$ is constant. Under suitable assumptions, we prove that the law of $u(t,x)$ has a density with respect to the
Lebesgue measure on $\R$. 

\begin{theorem}
\label{t6.1} We assume that $G$ satisfies the same assumptions than in Theorem \ref{t5.1} and that
\begin{equation*}
	\mathcal{J}(\delta):=\int_0^\delta  \Vert G(s,\ast)\Vert_\hac^2 ds>0,
\end{equation*}
for any $\delta>0$.
Suppose that the
function $\sigma$ is constant, and $b$ is continuously differentiable with bounded derivative. Then, for any $(t,x)\in]0,T]\times\Rd$ the probability law
of $u(t,x)$ has a density.
\end{theorem}
\noindent{\it Proof.} We apply Bouleau-Hirsch's criterion, see e.g. \cite[Section 2.1.3]{nualart}. 
 Fix $(t,x)\in ]0,T]\times \Rd$. We already know from Theorem \ref{t5.1} that
$u(t,x)\in\mathbb{D}^{1,2}$. Thus, it suffices to show that 
\begin{equation}
\label{6.c}
\Vert D u(t,x)\Vert_{\hac_T}^2 >0,\  a.s.
\end{equation}

From \eqref{5.1}, 
and for $\delta\in]0,t]$, we obtain
\begin{align}
\label{6.a}
\Vert Du(t,x)\Vert_{\hac_T}^2& = \int_0^t \Vert D_{s,\ast}u(t,x)\Vert_\hac^2 ds\ge \int_{t-\delta}^t \Vert D_{s,\ast}u(t,x)\Vert_\hac^2 ds\nonumber\\
&\ge\frac{1}{2}\sigma^2\int_{t-\delta}^t \Vert G(t-s,x-\ast)\Vert_\hac^2 ds
- I(t,x;\delta),
\end{align}
where
\begin{equation*}
I(t,x;\delta)= \int_{t-\delta}^\delta ds \left\Vert \int_0^t dr  \int_{\Rd} dz G(t-r,x-z) b^\prime(u(r,z))D_{s,\ast}u(r,z)\right\Vert_{\hac}^2.
\end{equation*}
By a change of variable, we have
\begin{equation}
\label{6.b}
\int_{t-\delta}^t \Vert G(t-s,x-\ast)\Vert_\hac^2 ds=\int_0^\delta  \Vert G(s,\ast)\Vert_\hac^2 ds = \mathcal{J}(\delta).
\end{equation}
Assumption {\bf(A3)} implies that $\mathcal{J}(T)<\infty$. Hence, 
\begin{equation}
\label{6.1}
\lim_{\delta\to 0} \mathcal{J}(\delta)=0.
\end{equation}
The Malliavin derivative $D_{s,\ast}u(r,z)$ vanishes except if $0\le s\le r$. Using this property and the change of variables $s\mapsto t-s$, $r\mapsto t-r$, we obtain
\begin{align*}
\E[I(t,x;\delta)]&=\E\left[\int_{t-\delta}^\delta ds \left\Vert \int_s^t dr  \int_{\Rd} dz G(t-r,x-z) b^\prime(u(r,z))D_{s,\ast}u(r,z)\right\Vert_{\hac}^2\right]\\
&= \E\left[\int_0^\delta ds \left\Vert \int_0^s dr  \int_{\Rd} dz G(r,x-z) b^\prime(u(t-r,z))D_{t-s,\ast}u(t-r,z)\right\Vert_{\hac}^2\right].
\end{align*}
We apply Fubini's theorem and then, \eqref{2.3} with $\mathcal{A}:=\hac$, $\Lambda:=G$ and $Z(r,z):=b^\prime(u(t-r,z))D_{t-s,\ast}u(t-r,z)$. Since the function $b^\prime$ is bounded, we obtain
\begin{align}
\label{6.2}
	\E[I(t,x;\delta)]	&	\le C \int_0^\delta ds\int_0^s dr \E\left[\left\Vert b^\prime (u(t-r,0) D_{t-s,\ast}u(t-r,0)\right\Vert_{\hac}^2\right] \sup_{\eta\in\Rd} \left\vert\tf G(\eta)\right\vert^2\nonumber\\
&\le C \E\left[\int_0^\delta ds \int_0^\delta dr \sup_{\eta\in\Rd} \left\vert\tf G(r)(\eta)\right\vert^2 \left\Vert D_{t-s,\ast} u(t-r,0)\right\Vert_\hac^2\right]\nonumber\\
&\le C  \int_0^\delta dr \sup_{\eta\in\Rd} \left\vert\tf G(r)(\eta)\right\vert^2 \E\left[ \left\Vert D_{t-\cdot,\ast} u(t-r,0)\right\Vert_{\hac_\delta}^2\right].
\end{align}
The next objective is to prove that
\begin{equation}
\label{6.3}
\sup_{0\le r\le \delta}\E\left[ \left\Vert D_{t-\cdot,\ast} u(t-r,0)\right\Vert_{\hac_\delta}^2\right]\le C \mathcal{J}(\delta).
\end{equation}
Indeed, owing to \eqref{5.1}, and by applying once more \eqref{2.3} as in \eqref{6.2}, we have
\begin{align*}
&\E \left[\left\Vert D_{t-\cdot,\ast} u(t-r,0)\right\Vert_{\hac_\delta}^2\right]\\
&\qquad \le 2\sigma^2 \mathcal{J}(\delta)
+ 2 \E\left[\left\Vert \int_0^{t-r}ds \int_{\Rd} dy G(t-r-s,z-y) b^\prime(u(s,y)) D_{t-\cdot,\ast}u(s,y)\right\Vert_{\hac_\delta}^2\right]\\
&\qquad\le C_1 \mathcal{J}(\delta) + C_2\int_0^{t-r} ds \left(\sup_{\eta\in\Rd}\left\vert \tf G(s)(\eta)\right\vert^2\right) \E\left[\left\Vert D_{t-\cdot,\ast}u(s,0)\right\Vert^2_{\hac_\delta}\right].
\end{align*}
Hence, \eqref{6.3} follows from an application of a version of Gronwall's Lemma.

From \eqref{6.2}, \eqref{6.3}, we obtain
\begin{equation}
\label{6.5}
\E[I(t,x;\delta)]\le C \mathcal{J}(\delta) \bar{\mathcal{J}}(\delta),%\int_0^\delta ds \sup_{\eta\in\Rd}\left\vert \tf G(s)(\eta)\right\vert^2.
\end{equation}
with 
\begin{equation*}
\bar{\mathcal{J}}(\delta):=\int_0^\delta ds \sup_{\eta\in\Rd}\left\vert \tf G(s)(\eta)\right\vert^2.
\end{equation*}
Notice that, assumption {\bf(A6)} on $G$ implies
\begin{equation*}
\lim_{\delta\to 0} \bar{\mathcal{J}}(\delta) = 0.
\end{equation*}

Fix $\delta\in]0,t[$ sufficiently small and $n\in\N$ sufficiently large such that $\frac{1}{n} <\frac{\sigma^2}{3} \mathcal{J}(\delta)$. Using Chebyshev's inequality along with \eqref{6.a}, \eqref{6.b}, \eqref{6.5} yield
\begin{align*}
\lim_{n\to\infty} \mathbb{P}\left[\Vert D u(t,x)\Vert_{\hac_T}^2 <\frac{1}{n}\right]&\le \lim_{n\to\infty}\mathbb{P}\left[ I(t,x;\delta)\ge \frac{\sigma^2}{2} \mathcal{J}(\delta)-\frac{1}{n}\right]\\
&\le \lim_{n\to\infty}\left(\frac{\sigma^2}{2} \mathcal{J}(\delta)-\frac{1}{n}\right)^{-1}\E\left[I(t,x;\delta)\right]\\
&\le C \bar{\mathcal{J}}(\delta).
\end{align*}
Letting $\delta\to 0$, we obtain
\begin{equation*}
\mathbb{P}\left[\Vert D u(t,x)\Vert_{\hac_T}^2 =0\right] = 0.
\end{equation*}
This is equivalent to \eqref{6.c}. 

\cqd
\medskip

Consider the particular case of the stochastic wave equation in spatial dimension $d>3$.
The Fourier transform of the fundamental solution of the corresponding partial differential equation is given by%We shall denote by $\Gamma$ the fundamental solution corresponding to the wave operator on $\Rd$. It is well known that
\begin{equation*}
\tf G(t)(\xi) = \frac{\sin(2\pi t |\xi|)}{2\pi |\xi|}.
\end{equation*}
Hence, there exist constants $C_1$, $C_2$, depending on $T$, such that for any $t\in[0,T]$, and $\xi\in\Rd$,
\begin{equation}
\label{6.6}
\frac{C_1}{1+|\xi|^2} \le \frac{\sin^2(2\pi t |\xi|)}{4\pi^2 |\xi|^2} \le \frac{C_2}{1+|\xi|^2}.
\end{equation} 
Assume that the spectral measure $\mu$ satisfies 
\begin{equation}
\label{6.7}
\sup_{\eta\in\Rd} \int_{\Rd}\frac{\mu(d\xi)}{1+|\xi + \eta|^2} < \infty.
\end{equation}
Then, according to \cite[Theorem 5.1]{conusdalang}, $G$ satisfies the conditions {\bf (A3)}, {\bf (A4)} and {\bf (A6)}. Hence, Theorem \ref{t5.1} holds.

Property \eqref{6.7} along with \eqref{6.6} imply
\begin{equation*}
\mathcal{J}(\delta)\ge C_1\int_0^\delta ds \int_{\Rd} \frac{\mu(d\xi)}{1+|\xi|^2}\ge C \delta.
\end{equation*}
This yields the following result
\begin{theorem}
\label{t6.2} Consider the particular case where $G$ is the fundamental solution of the wave equation with $d\in\N$. Assume \eqref{6.7} and that the functions $\sigma$ and $b$ are as in Theorem \ref{t6.1}.
Then, the statement of that theorem holds.
\end{theorem}

Assume that the covariance measure $\Gamma$ has a density: $\Gamma(dx)= f(x) dx$, with $f\ge 0$. In \cite{peszat}, it is proved that \eqref{6.7} is equivalent
to $\int_{\Rd} \frac{\mu(d\xi)}{1+|\xi|^2} < \infty$. This condition is satisfied for example if $f(x) =Ê|x|^{-\beta}$, $\beta\in]0,2[$, a case that has been extensively studied
in the literature of SPDEs driven by correlated noises.
\medskip

\noindent{\bf Remark} 
\smallskip

As has been already mentioned in the introduction, so far the existence of density for the probability law of the solution of an SPDE like \eqref{i.1} has been established when $G$
is a non-negative distribution. In this case, it is proved that the dominant term in the analysis of the Malliavin matrix is the first term on the right-hand side of \eqref{5.1}. Assuming that
the coefficient $|\sigma|\ge \sigma_0>0$, we have
\begin{equation}
\label{6.8}
\Vert G(t-\cdot,x-\ast) \sigma(u(\cdot,\ast))\Vert_{\hac_T}^2\ge \sigma_0^2 \Vert G(t-\cdot,x-\ast)\Vert_{\hac_T}^2.
\end{equation}
Then, the result is obtained by following a similar argument as in the proof of Theorem \ref{t6.1}.

For the wave operator in spatial dimension $d>3$, $G$ fails to satisfy the non-negativity requirement. So far, we have not been able to have a suitable lower bound like for instance in 
\eqref{6.8}. We notice that the trivial lower bound
\begin{align*}
\Vert G(t-\cdot,x-\ast) \sigma(u(\cdot,\ast))\Vert_{\hac_T}^2&=\int_0^t\int_{\Rd}|\tf G(t-s)(\xi)|^2 \mu_s^{\sigma(u)}(d\xi) ds\\
&\ge \sigma_0^2 \int_0^t \inf_{\eta\in\Rd} |\tf G(t-s)(\xi+\eta)|^2 \mu_s^{\sigma(u)}(d\xi) ds,
\end{align*}
does not help. Indeed, if $\mu$ is the spectral measure of a Riesz kernel ($\Gamma(dx)=|x|^{-\beta}dx$, $\beta\in]0,2[$), one can prove that the last integral in the above expression vanishes.

\bigskip

\noindent{\bf Acknowledgement} 
\smallskip

\noindent The authors thank the anonymous referees for their useful comments.

%%%%%End of the article

\end{document}